\documentclass[aos]{imsart}
\RequirePackage{natbib}
\usepackage{amssymb,amsmath,bm,graphicx,array,multirow,rotating,subfigure}
\usepackage{enumerate}
\usepackage{amsthm}

\newcounter{hypA}
\newenvironment{hypA}{\refstepcounter{hypA}\begin{itemize}
  \item[({A\arabic{hypA}})]}{\end{itemize}}
  \newcounter{hypB}

  \newcounter{hypC}

    \newcounter{hypD}

\def\cN{{\mathcal N}}

\def\bR{{\mathbb R}}

\def\1{{\mathbf 1}}
\def\eps{{\mathbf \varepsilon}}
\def\E{{\mathbf E}}

\def\P{{\mathbf P}}

\def\deq{{\stackrel{d}=}}

\def\cN{{\mathcal N}}
\def\cM{{\mathcal M}}
 
\newtheorem{theorem}{Theorem}
\newtheorem{lemma}{Lemma}

\newtheorem{remark}{Remark}

\begin{document}
\begin{frontmatter}

\title{High-dimensional change-point detection with sparse alternatives}
\runtitle{High-dimensional change-point detection}

\begin{aug}
\author{\fnms{Farida}  \snm{Enikeeva}\thanksref{t1,m1,m2}\ead[label=e1]{firstname.lastname@inria.fr}}
\and
\author{\fnms{Zaid} \snm{Harchaoui}\thanksref{t1,m1}\ead[label=e2]{firstname.lastname@inria.fr}}

 \thankstext{t1}{This work was partially funded by the Gargantua project under program Mastodons of CNRS, the LabEx PERSYVAL-Lab (ANR-11-LABX-0025) and by the OPTIMYST-II project of Minalogic.}

\runauthor{F. Enikeeva and Z. Harchaoui}

\affiliation{INRIA-LJK\thanksmark{m1} and IITP\thanksmark{m2}}
\address{INRIA Grenoble Rh\^one-Alpes\\ 655, Avenue de l'Europe, Montbonnot\\ 38334 Saint-Ismier Cedex, France \\ 
          \printead{e1,e2}}
\address{Institute for Information Transmission Problems\\ Bolshoy Karetny per. 19, Moscow, 127994, Russia}
\end{aug}

\begin{abstract}
We consider the problem of detecting a change in mean in a sequence of Gaussian vectors. Under the alternative hypothesis, the change occurs only in some subset of the components of the vector. We propose a test of the presence of a change-point that is adaptive to the number of changing components. Under the assumption that the vector dimension tends to infinity and the length of the sequence grows slower than the dimension of the signal, we obtain the detection boundary for this problem and prove its rate-optimality.
\end{abstract}

\begin{keyword}[class=AMS]
\kwd[Primary ]{62G10}
\kwd{62H15}
\kwd[; secondary ]{60G35}
\end{keyword}

\begin{keyword}
\kwd{change-point problem, high-dimensional data, sparsity, detection boundary}
\end{keyword}

\end{frontmatter}

\section{Introduction}
Consider a sequence of independent $d$-dimensional Gaussian vectors $X_1,\dots,X_n$ with a possible change in mean at an unknown location $\tau$,
\begin{equation}\label{channels}
X_i=\theta+\Delta\theta_\tau\1\{ i>\tau \}+\xi_i,\quad i=1,\dots,n.
\end{equation}
Here $\xi_i$ are i.i.d. random variables drawn from $\cN(0,I_d)$, $\Delta\theta_\tau\in\bR^d$, $\theta\in\bR^d$. Our goal is to design a test for the change-point problem
\begin{align*}
\mathbf{H}_{0}: &\quad \tau=n\,,\quad \text{and }\Delta\theta_\tau=\mathbf{0} \\
\mathbf{H}_{A}: &\quad \tau\in\{1,\dots,n-1\},\quad \text{and }\Delta\theta_\tau\neq\mathbf{0} \text{ with } 1\le p\le d,
\end{align*}
where $p$ stands for the number of non-zero components of the vector $\Delta\theta_\tau$. Under the null hypothesis $\mathbf{H}_{0}$, there is no change in mean ($\tau=n$) and~$\Delta\theta_\tau=\mathbf{0}$. Under the alternative hypothesis $\mathbf{H}_{A}$, a change occurs at the moment $\tau$ on a subset of components of the mean vector corresponding to $\mathrm{supp}(\Delta\theta_\tau)\subseteq\{1,\dots,d\}$. The alternative is composite, since neither the set $\text{supp}(\Delta\theta_\tau)$ nor its dimension $p$ are known. 

Change-point problems for multivariate Gaussian observations have received a lot of attention for decades. A classical approach consists in assuming that the change-point occurs on all components, \text{i.e.}  $\text{supp}(\Delta\theta_\tau) = \{1,\dots,d\}$, and to study the problem under a classical asymptotic setting, that is by letting the number of observations grow to infinity $n\to\infty$ while the dimension
of the vector $d$ remains fixed. We refer the reader to, \text{e.g.}~\citep{Shiryaev:1978, Basseville:Nikiforov:1993, Carlstein:Mueller:Siegmund:1994, Brodsky&Darhovsky:1993, Csorgo&Horvath:1997, Lehmann:Romano:2005, korostelev:2011} for a review. 

In this paper, we consider high-dimensional observations, and work in a double-asymptotic framework where both the number of observations and their dimension grow to infinity. Such a setting is particularly appropriate for recent instances of change-problems that arise in real-world applications, 
for example, in bioinformatics (copy-number variation analysis,~\cite{Bleakley&Vert:2011,Zhang&etal:2010}), in network traffic data analysis~\citep{LR2009,LLC:2012}, in multimedia indexation~\citep{HarchaouiVLC09},
and in finance. In all these applications, one is interested in detecting change-points in relatively short sequences of observations (say, $n=100$) whereas the  dimension of these observations is high (say, $d=10^3$). However, in these  applications, prior information suggests that the change is likely to occur \emph{only in a small subset of components}. Therefore,  the effective dimension  of
the change-point problem is actually the dimension of $\text{supp}(\Delta\theta_\tau)$, instead  of $d$, that makes the statistical problem potentially tractable even for very high-dimensional observations. 

\cite{Korostelev&Lepski:2008} studied high-dimensional change-point problems in the white noise framework. Yet, the authors assume that a change occurs in all components, \textit{i.e}~$p=d$. \cite{Korostelev&Lepski:2008}  propose an asymptotically minimax estimator of the
change-point location, under double-asymptotics with $\|\Delta\theta\|\to\infty$ as $d\to\infty$.  See also~\cite{opac-b1092742} for an early treatment of a related problem. Recently,~\cite{Xie&Siegmund:2013} consider a  problem similar to ours from a methodological point of view.
In~\citep{Xie&Siegmund:2013}, a Bayes-type test statistic is proposed, where the authors introduce a mixture model that hypothesizes an assumed fraction of changing components.

In this work, we compute the detection boundary for this change-point problem. The detection boundary is an asymptotic condition on the norm of $\Delta\theta_\tau$ providing the minimax separability of the hypotheses ${\mathbf H}_0$ and $\mathbf{H}_A$. It depends on the location  of the change-point and on the number of changing components $p$. The proposed test is based on two base test statistics: a linear  statistic, that
considers all components simultaneously; a scan statistic that searches for a change over all possible combinations of changing components. Although the problem has a combinatorial structure, we show that the proposed test statistic can be computed efficiently, in almost linear time with respect to the dimensions of the problem. We derive the rate of convergence of our test and prove that it is adaptive to the unknown set of changing components and that it is rate-optimal.

\section{Statement of the problem}
We first consider the problem of testing the hypothesis of no-change against the alternative of a change in mean at location $\tau$ in exactly $p$ components. We will describe later the test that is adaptive to the case of unknown $p$ and $\tau$. 

Both the null and the alternative hypothesis can be simply formulated in terms of the norm of jumps of the mean vector, $\Delta\theta_\tau$, so that we have no change if the norm of jumps is zero, $\|\Delta\theta_\tau\|=0$. We say that the change occurs at location $\tau$, if the norm of  jumps at location $\tau$ satisfies $\|\Delta\theta_\tau\|>r$ for some value $r>0$. 

Define the set $\Theta_{p}[r]=\{v\in V_p^d:\ \|v\|\ge r \}$, and 
$$
V_p^d=\{v=(\eps_1v_1,\dots,\eps_dv_d):\ \eps_j\in\{0,1\},\ \sum\limits_{j=1}^d \eps_j=p\}
$$ 
a subspace of $\bR^d$-vectors with $p$ non-zero components, with $\|\cdot\|$ denotes the Euclidean norm. 

Let $I\subseteq \{1,\dots,n-1\}$ be a set of possible locations of a change. We consider two general problems based on model~(\ref{channels}). 
\begin{enumerate}[(P1)]
\item Testing a change in mean in exactly $p$ components at a given location $\tau\in I$:
$$
\mathbf{H}_{0}:\ \Delta\theta_\tau=0\quad \mbox{vs}\quad {\mathbf{H}_{A}}:\ \Delta\theta_\tau \in \Theta_{p}[r],
$$
where $r>0$ may depend on $\tau$, $p$, $d$, and $n$.
\item Testing the presence of a change-point in an unknown number of components within the time interval $I$:
$$
\mathbf{H}_{0}:\ \Delta\theta_\tau=0\ \forall \tau\in \{1,\dots,n-1\} \quad \mbox{vs}\quad {\mathbf{H}_{A}}:\ \exists \tau\in I:\ \Delta\theta_\tau \in \bigcup\limits_{p=1}^d \Theta_p[r].
$$
Note that $r$ might depend on $\tau$, $p$, $d$ et $n$ as well; 
we shall however omit to explicit this dependence to keep light notations. 
\end{enumerate} 

Note that Problem~(P1) corresponds to a two-sample test problem with mean-difference in a subset of components and same variance; 
see~\citep{Tony_Cai_Liu_Xia_2014} for a review of recent works on this topic.

A test $\psi=\psi(X_1,\dots,X_n)$ is a measurable function of observations~(\ref{channels}).  
For any test $\psi$, define the type I error as $\alpha(\psi)= \E_0 \psi$.  Type II errors  for Problems~(P1) and~(P2) are respectively defined as
\begin{align*}
\beta(\psi,\Theta_p[r],\tau)&=\sup_{\Delta\theta_\tau\in\Theta_p[r]} \E_{\Delta\theta_\tau} (1-\psi ) \\
\beta^{*}(\psi,\Theta[r], I)&=  \sup_{\tau\in I} \sup_{\Delta\theta_\tau\in\Theta[r]} \E_{\Delta\theta_\tau}(1-\psi )\\
&\equiv \sup_{\tau\in I} \sup_{p=1,\dots,d}\beta(\psi,\Theta_p[r],\tau).
\end{align*}
where $\P_{\Delta\theta_\tau}$ is the distribution corresponding to the alternative with a $\Delta\theta_\tau$ change in mean at location $\tau$.  
Define the global testing errors~\citep{Ingster&Suslina:2003} for two problems:
\begin{align*}
\gamma(\psi,\Theta_p[r],\tau)&:= \alpha(\psi)+\beta(\psi,\Theta_p[r],\tau) \\
\gamma^{*}(\psi,\Theta[r],I)&:= \alpha(\psi)+\beta^{*}(\psi,\Theta[r],I),
\end{align*}
where $\Theta[r]:=\bigcup\limits_{p=1}^d \Theta_{p}[r]$.

We make the following Assumptions (A1-A2-A3) on the asymptotic behaviour of $n$, $d$, $p$, and $\tau$. 
\textit{Throughout the paper, the asymptotics of $p$ and $n$ are parametrized by $d$, where $d\to\infty$.}
The asymptotic of the location of $\tau$ is naturally parametrized by $n$, which, in turn, depends on $d$. 
\begin{hypA}\label{asymp_gal}
The number of observations grows with the number of components
\begin{equation*}
 n=n(d)\to\infty\, \text{ as }d\to\infty.
\end{equation*}
\end{hypA}
\begin{hypA}\label{asymp_cond1}
The number of components with a change is sufficiently large
\begin{equation*}
p\to \infty\quad \mbox{and}\quad d/p\to\infty\quad \mbox{as $d\to\infty$} .
\end{equation*} 
\end{hypA}
\begin{hypA}\label{asymp_cond2}
For Problem~(P2), we need an additional assumption, namely that
\begin{equation*}
\lim_{d\to\infty} \frac{\log (nd)}{p\log (d/p)}=0.
\end{equation*}
In particular, this assumption  implies that $\log (nd)/ \log {d\choose p}\to 0$ and $\log n/d\to 0$ as $d\to\infty$. 
\end{hypA}

\begin{remark}
If $p$ depends on $d$ via a sparsity coefficient $\beta\in [0,1)$, $p\asymp d^{1-\beta}$, then Assumption~(A\ref{asymp_cond1}) is satisfied. We shall distinguish between the cases of {\it high sparsity}, $\beta\in(1/2,1)$ and {\it low sparsity}, $\beta\in[0,1/2]$. 

\end{remark}

We are interested in the minimax separation conditions for problems~(P1) and~(P2).  The question is how far from the origin should be the sets $\Theta_p[r]$ and $\Theta[r]$  in order to separate the hypotheses $\mathbf{H}_{0}$ and ${\mathbf{H}_{A}}$ in problems~(P1) and~(P2), respectively. The sequence $r=r(d)\equiv C r_d$ is called {\it a detection boundary}~\citep{Ingster:1997} for problem~(P1) if 
\begin{enumerate}[(i)]
\item for any small $\alpha>0$ there exists a constant $C^*$ and a decision rule $\psi^*$ such that  for $r= C^* r_d$
$$
\limsup_{d\to\infty} \gamma (\psi^*,\Theta_p[r],\tau) \le \alpha
$$
\item there exist positive constants $C_*$ and $r_*$ such that for $r\equiv C_* r_d$ and any test $\psi$
$$
\liminf_{d\to\infty} \gamma(\psi,\Theta_p[r],\tau) \ge r_*
$$
\end{enumerate}
We can similarly define the detection boundary for problem~({P2}). The quantity $r_d$ is called {\it the minimax separation rate};
we refer to \citep{Ingster&Suslina:2003, korostelev:2011, Baraud:2002, donoho2004} for further discussion about this definition.

We are interested in the asymptotic conditions on the detection boundary $r= C r_d$ that
provide the minimax separation between the hypotheses as well as in the detection boundary constant $C$. 
The detection boundary for Problems~(P1) and~(P2) depends on the location of the change-point via the function
\begin{equation}\label{htau}
h(\tau)=\frac \tau n \Bigl(1-\frac\tau n\Bigr) \, ,
\end{equation}
which usually appears in change-point problems; see~\citep{Csorgo&Horvath:1997} for details.


We  use the following notation throughout the text:  $\cM(d,p)$ stands for the collection of all subsets of $\{1,\dots , d\}$ of cardinality $p$, 
and  $\cM$ stands for the set of all possible subsets of $\{1,\dots,d\}$. Denote by $\Pi_m v$, the projection of a vector $v\in\bR^d$ onto a subspace indexed by $m\in\cM$, and by $\varkappa$ the constant for which $\varkappa>2/(1-\log 2)$ (the latter constant, 
which appears in all results, is related to deviations of quadratic forms, and first appears in~Lemma~\ref{Spok}).

\section{Testing procedure}
Let us first define the $d$-dimensional process $Z_n(s)$, $s=1,\dots,n-1$, where
\begin{equation}\label{Zj}
Z_n(s)=\sqrt{\frac{s(n-s)}{n}}\biggl(\frac1s\sum_{i=1}^s X_i-\frac1{n-s}\sum_{i=s+1}^n X_i\biggr).
\end{equation}
The proposed tests for problems (P1)--(P2) are based on two $\chi^2$-type test statistics, which will be referred to as \emph{linear statistic} and \emph{scan statistic}. 

Let $s\in \{1,\dots,n-1\}$. The linear statistic is given by
\begin{equation}\label{lin_stat}
L_{\mathrm {lin}}(s)=\frac{\|Z_n(s)\|^2-d}{\sqrt{2d}}.
\end{equation}

For each fixed $p\in\{1,\dots,d\}$ the scan statistic is defined as
\begin{equation}\label{scan_stat_p}
L^p_{\mathrm {scan}}(s)=\max_{m\in\cM(d,p)} \left\{\frac{\|\Pi_m Z_n(s)\|^2-p}{\sqrt{2p}}\right\}.
\end{equation}

\begin{remark}
Assume that $\Delta\theta_\tau\in V_p^d$ and let $m\in\cM(d,p)$ be a subset of $p$ components with a change. The choice of the test statistics based on $Z_n(s)$ is motivated by the fact that the marginal log-likelihood ratio of  $\tau\in\{1,\dots,n-1\}$ and $m$ is given by
\begin{equation}\label{LLZ}
\log\frac{d\P_{\tau,m}}{d\P_0}(X)=\frac12\|\Pi_m Z_n(\tau)\|^2,
\end{equation}
where $\P_{\tau,m}$ is the measure corresponding to the presence of a change at the location $\tau$ in the set of components $m$. 
The proof can be found in Lemma~\ref{LL_lemma} in Appendix~C.
\end{remark}

\subsection{Known number of components with a change and known $\tau$}
The proposed decision rule for Problem~(P1) is based on the combination of two tests 
$$
\psi_p^*=\psi_{\mathrm {lin}}\vee\psi^p_{\mathrm {scan}},
$$
with  
$$
\psi_{\mathrm {lin}}=\1\bigl\{L_{\mathrm {lin}}(\tau)>H\bigr\},\quad 
\psi^p_{\mathrm {scan}}=\1\bigl\{L^p_{\mathrm {scan}}(\tau )>T_p\bigr\}.
$$
Thresholds $H$ and $T_p$ should be set so that the global risk error $\gamma(\psi,\Theta_p[r])$ tends to zero
as $d \to \infty$. Theorem~\ref{up_thm} answers this question, and provides a guideline to set the corresponding thresholds
depending on $d$ (and henceforth on $p$ and $n$). 
Should we work in a Neyman-Pearson setting, for any significance level $\alpha$ we can set the thresholds for the two tests 
using the quantiles of $\chi^2$ distributions
$$
H=(\chi_d^2)^{-1}\Bigl[1-\alpha \Bigr],  \quad T_p=(\chi_p^2)^{-1}\left[1-\alpha\Bigl/{d\choose p}\right].
$$
However, such a strategy imposes a computational burden for the scan test, since the quantiles of such a high order could be difficult to compute precisely even for small values of $d$. We propose the following formulas for the thresholds
\begin{equation}\label{approx_thresh_formulas}
H=\sqrt{(\varkappa/2)\log (1/\alpha)},\quad T_p=\frac \varkappa {\sqrt{2 p}}\log \biggl[{d\choose p} \frac{1}\alpha\biggr]
\end{equation}
using concentration inequalities for the norm of $d$-dimensional Gaussian vector (see Remark~\ref{thrA} to Lemma~\ref{type_Ierr} for the details about these thresholds).

\begin{remark}
The scan statistic $\psi^p_{\mathrm {scan}}$ and the linear test statistic $\psi_{\mathrm {lin}}$ have different performance depending on the sparsity index. In the case of high sparsity, $\beta\in(1/2,1)$,  the scan statistic test outperforms the linear test, since the scan statistic searches for a change over all possible combinations.  However, in the case of moderate sparsity, $\beta\in [0,1/2]$, 
the linear has faster detection rate. In other words, pooling statistical information from all vector components is more effective than a full search over all possible combinations of components. Therefore, the proposed test statistic $\psi_p^*$ gets the best of both worlds, and is effective in all situations, moderate sparsity and high sparsity. 
\end{remark}

\subsection{Adaptation to unknown number of components with a change}
For Problem~(P2), the proposed adaptive decision rule is, again, based on the combination of two tests 
$$
\psi^*=\psi_{\mathrm {lin}}\vee\psi_{\mathrm {scan}}.
$$
Here, note that the linear test statistic does not change. However, for Problem (P2), the scan statistic is defined as 
\begin{equation}\label{scan_stat}
L_{\mathrm {scan}}(s)=\max_{p=1,\dots,d} \frac1{T_p} \max_{m\in\cM(d,p)} \left\{\frac{\|\Pi_m Z_n(s)\|^2-p}{\sqrt{2p}}\right\}
\end{equation}
where 
\begin{equation}\label{scan_threshold}
T_p=\frac \varkappa {\sqrt{2 p}}\log \biggl[{d\choose p} \frac{nd}\alpha\biggr].
\end{equation}
Again, for any significance level $\alpha$, we can derive for the two tests
$$
\psi_{\mathrm {lin}}=\1\left\{\max_{s\in I} L_{\mathrm {lin}}(s)>H\right\},\quad 
\psi_{\mathrm {scan}}=\1\left\{\max_{s\in I} L_{\mathrm {scan}}(s)>1\right\}
$$
approximate formulas similar to the ones in~(\ref{approx_thresh_formulas});
see Remark~\ref{rem:approx_thresh} in Appendix~A for details.
Note that here the threshold $T_p$ is built in the scan test statistic. 

\begin{remark}
At first sight, the scan statistic seem to be difficult to compute, since it involves a combinatorial search which could be computationally hard. Recent work~\citep{addario-berry2010} showed several general classes of problems
where scan statistics are computationally hard to compute. However, in our problem, it turns out that the scan statistic can be efficiently computed in almost linear-time with respect to the dimension of the problem. Indeed, we have
$$
\max_{m\in\cM(d,p)} \frac{\|\Pi_m Z_n(s)\|^2-p}{\sqrt{2p}} =\frac1{\sqrt{2p}}\biggl(\sum_{j=1}^p [Z_n^{(j)}(s)]^2 -p\biggr),
$$
where $[Z_n^{(j)}(s)]^2$ are the ordered squred components of the vector $Z_n(s)]$:
$$
[Z_n^{(1)}(s)]^2>[Z_n^{(2)}(s)]^2>\dots >[Z_n^{(d)}(s)]^2.
$$
Thus the computational complexity of the adaptive test statistic is $O(nd\log d)$. Computation rely on sorting the squared components of vectors $Z_n(s)$ for each $s=1,\dots,n-1$. 
\end{remark}

\begin{remark}
Again, the proposed adaptive test statistic $\psi^*$ covers all situations, both moderate sparsity and high sparsity. This is reflected by our theoretical results in Theorems~\ref{up_thm_adapt}-\ref{low_thm_adapt}, which show that the proposed adaptive test statistic is minimax-optimal and rate-adaptive. 
Simulations in Section~\ref{sec:exps} corroborate our theoretical results. 
\end{remark}

\section{Main results}
\label{sec:mainresults}

\subsection{Upper and lower bounds in problem (P1)} 
The following theorem gives the upper bound for the test $\psi_p^*$ in problem~(P1). 
\begin{theorem}\label{up_thm}
Assume that Assumptions~(A\ref{asymp_gal})--(A\ref{asymp_cond1}) hold and $r=r(d)$ satisfies either
\begin{equation}\label{lin_dbp}
\lim_{d\to\infty} \frac{r^2 n h(\tau)}{\sqrt d}=+\infty
\end{equation}
or
\begin{equation}\label{scan_dbp}
\liminf_{d\to\infty} \frac{r^2 nh(\tau)}{p\log (d/p)}>\varkappa\ge 6.6.
\end{equation}
Let $H\to\infty$  as $d\to \infty$ such that $H\ge\sqrt{(\varkappa/2)\log d}$ and 
\begin{equation}\label{H}
\limsup_{d\to\infty} \frac{H \sqrt{2d}}{r^2 nh(\tau)} <\sqrt 2-1.
\end{equation}
Let for some $\delta>0$
\begin{equation}\label{Tp}
T_p=\frac \varkappa {\sqrt{2 p}}\log \biggl[(1+\delta) {d\choose p}\biggr].
\end{equation}
Then $\gamma(\psi_p^*,\Theta_p[r],\tau)\to 0$  as $d\to\infty$.
\end{theorem}

The following theorem establishes the lower bound on the minimax error probability in problem~(P1) 
and the minimax separation rates providing the separability conditions between the hypotheses. We also obtain the detection boundary constant.
\begin{theorem}\label{low_thm}
Assume that $p\asymp  d^{1-\beta}$ as $d\to\infty$, Assumption~(A\ref{asymp_gal}) holds and $r=r(d)$ satisfies two following conditions
\begin{align}
&\lim_{d\to\infty} \frac{r^2 n h(\tau)}{\sqrt d}=0,\label{det_bound_lin}\\
&\limsup_{d\to\infty} \frac{r^2 nh(\tau)}{p\log (d/p)}<2-\frac1\beta,\quad \beta\in (1/2,1).\label{det_bound_scan}
\end{align}
Then  $\inf\limits_{\psi\in[0,1]}\gamma(\psi,\Theta_p[r],\tau)\to 1$ as $d\to\infty$.
\end{theorem}

\begin{remark}
Consider the case of a single vector, $n=1$. Then, our problem could be rephrased as testing $p$ non-zero components in the mean of a Gaussian vector. Such a problem was previously considered by~\citep{Baraud:2002} in non-asymptotic setting and by \citep{Ingster:2002b}. In this particular problem, we recover in~(\ref{det_bound_lin})--(\ref{det_bound_scan}), a detection boundary which is similar to the one given in~\citep{Baraud:2002}. The quantity $2-1/\beta$ in~(\ref{det_bound_scan}) coincide with the key quantities arising in the problem of classification of a Gaussian vector with $p=d^{1-\beta}$ non-zero components in the mean obtained in~\citep{Ingster&Pouet&Tsybakov:2009} and in the problem of detection~\citep{Ingster:2002b, Ingster:1997}. Recently,~\cite{butucea2013} obtained similar results for the problem of detection of a sparse submatrix of a noisy matrix of growing size.
\end{remark}

\subsection{Adaptation}

The following theorem gives the upper bound for the adaptive test $\psi^*=\psi_{\mathrm{lin}}\vee \psi_{\mathrm{scan}}$. 
\begin{theorem}\label{up_thm_adapt}
Assume that (A\ref{asymp_gal}--A\ref{asymp_cond2}) hold and $r=r(d)$ satisfies 
\begin{equation}\label{up_adapt_lin}
\lim_{d\to\infty} \min_{\tau\in I} \frac{r^2n h(\tau)}{\sqrt{2d\log(d\log n)}} =+\infty
\end{equation}
or
\begin{equation}\label{up_adapt_scan}
\liminf_{d\to\infty} \min_{\tau=1,\dots,n-1}\min_{p=1,\dots,d}\frac{r^2 nh(\tau)}{p\log (d/p)}>\varkappa\ge 6.6.
\end{equation}
Let $c_0>\sqrt{\varkappa}$ and $c_1\in(0,\sqrt 2-1)$. Choose the thresholds for the linear test as
\begin{equation}\label{H_adapt}
c_0\sqrt{ \log (d\log n)} \le H\le (\sqrt 2-1-c_1) \min_{s\in I } \frac{r^2 nh(s)}{\sqrt{2d}}
\end{equation}
and for the scan test as given by
$$
T_p=\frac{\varkappa}{\sqrt{2p}} \log\left[{d\choose p} |I| dp^2\right].
$$
Then $\gamma^*(\psi^*,\Theta[r],I)\to 0$ as $d\to\infty$.
\end{theorem}
The lower bound is a direct consequence of Theorem~\ref{low_thm}. Recall that $\Theta[r]=\bigcup\limits_{p=1} ^d\bigcup\limits_{\tau\in I} \Theta_p[r]$.
\begin{theorem}\label{low_thm_adapt}
Assume that (A\ref{asymp_gal}) holds.. If there exist $\tau_0\in I$ and  $p_0\in\{1,\dots,d\}$, $p_0\asymp d^{1-\beta_0}$ as $d\to\infty$  such that 
$$
\lim_{d\to\infty} \frac{r^2 n h(\tau_0)}{\sqrt d}=0
$$
and 
$$
\limsup_{d\to\infty} \frac{r^2 nh(\tau_0)}{p_0\log (d/p_0)}<2-\frac1\beta_0,\quad \beta_0\in(1/2,1)
$$
then $\inf\limits_{\psi\in[0,1]}\gamma^*(\psi,\Theta[r],I)\to 1$ as $d\to\infty$.
\end{theorem}

\begin{remark}
Assuming that condition (A\ref{asymp_cond2}) is satisfied, Theorems~\ref{up_thm_adapt}-\ref{low_thm_adapt} show that
the proposed scan test reaches the minimax-optimal rate of convergence with no loss in the detection boundary rate. 
Regarding the proposed linear test, we get only a $\log(d\log n)$ loss in the detection boundary rate.  
It would be interesting to investigate counterparts of our approach based on higher-criticism ideas~\cite{donoho2004}, 
to see if it could bridge the gap between the constants arising in the upper and lower bounds. 
\end{remark}

\subsection{Discussion: detection boundary and rates}
Suppose that the change size in each component is constant, $\Delta\theta_j=a$, $j=1,\dots,d$. Thus the squared norm of the jumps of the vector mean $\|\Delta\theta_\tau\|^2$ equals $p a^2$ if we have $p$ components with a change. The detection boundary conditions can now  be written in terms of asymptotic behaviour of the jump size $a_{n,d}$. Recall that $p\asymp d^{1-\beta}$ where $\beta\in [0,1)$ is the sparsity coefficient. The rates depend on the location of the change point via the function $h(\tau)=\tau/n(1-\tau/n)$. 

As follows from Theorems~\ref{up_thm} and~\ref{low_thm}, in the case of moderate sparsity, $\beta\in [ 0,1/2]$, the detection boundary is of the form 
$$
 a_{d,n}^2 d^{1/2-\beta} nh(\tau) \asymp 1.
$$
On the other hand, for the case of high sparsity, $\beta\in(1/2,1)$, the detection is impossible if 
$$
\limsup_{d\to\infty} a_{d,n}^2 \frac{nh(\tau)}{\log d} <2\beta-1.
$$
We can express the detection boundary condition in the following way. Suppose that 
$$
a_{d,n}=r_d\sqrt{\frac{\log d}{nh(\tau)}}.
$$
 Then the detection is impossible if $\limsup\limits_{d\to\infty} r_d<\sqrt{2\beta-1}$  
 as shown in Theorem~\ref{low_thm}. On the other hand, from Theorem~\ref{up_thm} it follows that the detection is always possible if $\liminf\limits_{d\to\infty} r_d>\sqrt{\varkappa\beta}$. We see that there is a gap between the constants in the detection boundary. The main results on the detection boundaries are gathered in the table below. 
 
 \footnotesize{
\begin{center}
\begin{tabular}{|c|c|c|}
\hline \rule[-2ex]{0pt}{5.5ex}  & high sparsity, $\beta\in(1/2,1)$  &   moderate sparsity, $\beta\in [0,1/2]$\\ 
\hline \rule[-2ex]{0pt}{5.5ex} {\bf detection boundary} & $a_{d,n}=r_d\sqrt{\frac{\log d}{nh(\tau)}}$  &  $a_{d,n}^2 d^{1/2-\beta} nh(\tau) \asymp 1$\\ 
\hline \rule[-2ex]{0pt}{5.5ex}  upper bound & $\liminf\limits_{d\to\infty} r_d>\sqrt{\varkappa\beta}$ & $\lim\limits_{d\to\infty} a_{d,n}^2 d^{1/2-\beta} nh(\tau) =+\infty$\\ 
\hline \rule[-2ex]{0pt}{5.5ex}  lower bound & $\limsup\limits_{d\to\infty} r_d<\sqrt{2\beta-1}$  &  $\lim\limits_{d\to\infty} a_{d,n}^2 d^{1/2-\beta} nh(\tau) =0$\\ 
\hline \rule[-2ex]{0pt}{5.5ex} adaptive upper bound &  $a_{d,n}=r_d\sqrt{\frac{\log d}{nh(\tau)}}$ & $a_{d,n}^2 \frac{d^{1/2-\beta} nh(\tau)}{\sqrt{\log (d\log n)}} \asymp 1$ \\ 
\hline 
\end{tabular} 
\end{center}
}
\medskip
\normalsize{
We can consider three possible frameworks that take into account the location of the change-point:
\begin{enumerate}[(i)]
\item $n$ is fixed, $d\to\infty$, the change occurs at $\tau\in\{1,\dots,n-1\}$ and $\tau/n=t^*\in(0,1)$;
\item $n=n(d)$ grows as $d\to\infty$ and the change occurs at $\tau$ such that $\tau/n\to t^*\in(0,1)$ and  $h(\tau)\to t^*(1-t^*)\in(0,1)$ as $n\to\infty$
\item $n=n(d)$ grows as $d\to\infty$ and the change is close to the beginning or the end of the sample: $h(\tau) \asymp n^{-\delta}$ as $n\to\infty$, where $\delta>0$.
\end{enumerate}
The following table provides the asymptotics of the detection boundaries for~frameworks~(i)--(iii). Here $r_d'=r_d/\sqrt{nt^*(1-t^*)}$ is the renormalized detection boundary constant for framework~(i).
}
\footnotesize{
\begin{center}
\begin{tabular}{|c|c|c|}
\hline \rule[-2ex]{0pt}{5.5ex}  & high sparsity, $\beta\in(1/2,1)$  &  moderate sparsity, $\beta\in [0,1/2]$\\ 
\hline \rule[-2ex]{0pt}{5.5ex}  framework (i) &  $a_{d} = r_d'\sqrt{\log d}$ & $ a_d^2 d^{1/2-\beta} \asymp 1$ \\ 
\hline \rule[-2ex]{0pt}{5.5ex} framework~(ii) & $a_{d,n}=r_d\sqrt{\frac{\log d}{n t^*(1-t^*)}}$  &  $a_{d,n}^2 d^{1/2-\beta} n \asymp 1$\\ 
\hline \rule[-2ex]{0pt}{5.5ex}   framework~(iii) & $a_{d,n}=r_d \sqrt{\frac{\log d}{n^{1-\delta }}}$& $a_{d,n}^2 d^{1/2-\beta} n^{1-\delta} \asymp 1$\\ 
\hline 
\end{tabular} 
\end{center}
}
\normalsize{}
\section{Simulations}
\label{sec:exps}
In this section, we perform simulations to evaluate the empirical behaviour of the proposed test statistics, 
depending on the number of obserations $n$, the sparsity index $p$, the renormalized size of the jumps $\|\Delta\theta\|/\sqrt p$, and the dimension $d$. 
We also propose and evaluate two stategies for calibrating the proposed test statistics. 

First, we generated $n=100$ $d$-dimensional Gaussian independent vectors with independent components and a change in mean in $p$ components. We considered several situations of the change in mean in all components, $p=d$, the case of moderate sparsity, $p\gg \sqrt d$, and the case of high sparsity, $p\ll \sqrt d$. We performed $500$ repetitions, and report results averaged over these repetitions. The significance level was set to $\alpha=0.05$.  In all simulations except the last one, we used the $\alpha/2$ empirical quantile of the $\chi^2(d)$-distribution instead of $H$ for the calibration of the linear test at level $\alpha/2$. We used the $\alpha/(2d)$ empirical quantiles of the $\chi^2(p)$-distribution to calibrate the scan statistic test at level $\alpha/2$.

\begin{figure}[htbp!]
\includegraphics[width=0.8\textwidth]{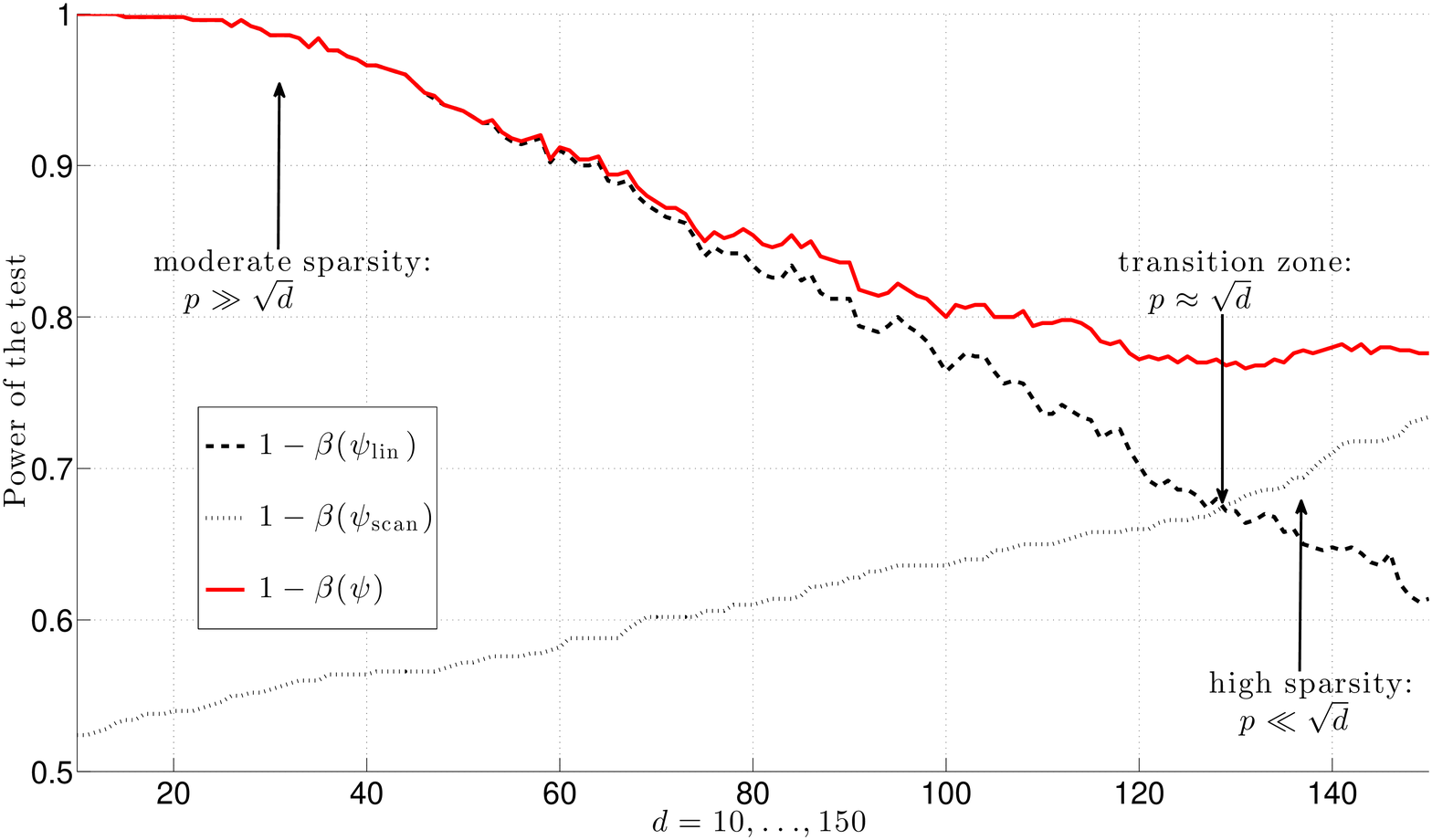}
\caption{Power of three tests for $n=100$, $p=10$, $\tau=n/4$, $\Delta\theta_j=0.6$.}
\label{fig:1}
\end{figure}

In Figure~\ref{fig:1} the power of the test is shown for the linear test $\psi_{\mathrm {lin}}$, 
the scan test $\psi_{\mathrm {scan}}$ and the proposed test $\psi$. The number of channels $d$ varies from 10 to 150 while $p=10$ is fixed. The mean in each channels changes to $\Delta\theta_j=\pm 0.6$ at the location $\tau=n/4$. We can  that the scan statistic performs better in the case of high sparsity whereas the power of the linear statistic test is much higher in the case of the moderate sparsity and in the non-sparse case. At the same time in the transition zone $p\approx \sqrt{d}$ both test have a similar performance.

\begin{figure}[htbp!]
 \includegraphics[width=0.8\textwidth]{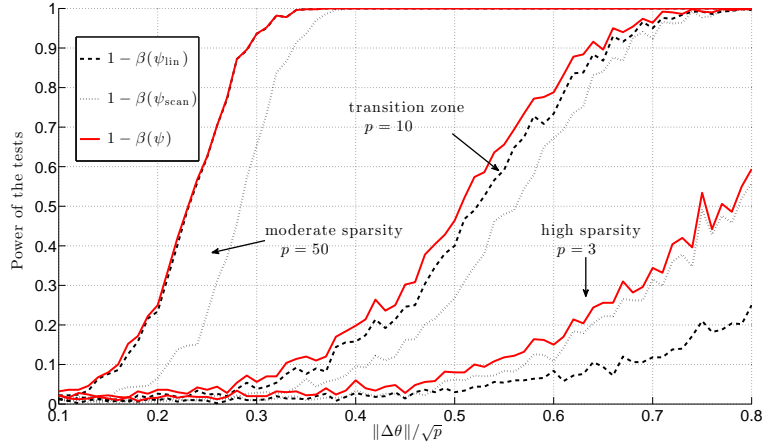}
 \caption{Power of the tests for $n=100$, $d=100$, $\tau=n/4$ for different values of $p$.}
 \label{fig:2}
 \end{figure}

We also ran simulations for $d=100$, with several sparsity indices: $p=3,\ 10,\ 50$. In Figure~\ref{fig:2},
we report the empirical power of the three tests depending on $\|\Delta\theta\|/\sqrt p$. 
We observe that the detection boundary constant, which is proportional to $\|\Delta\theta\|$, increases as $p$ decreases. 
In the high sparsity case ($p=3$), we observe that the test based on the scan statistics outperforms the linear test. 
On the other hand, in the moderate sparsity case ($p=50$), the linear test works better than the scan statistic test.
 
\begin{figure}[htbp!]
\centering
\includegraphics[width=0.8\textwidth]{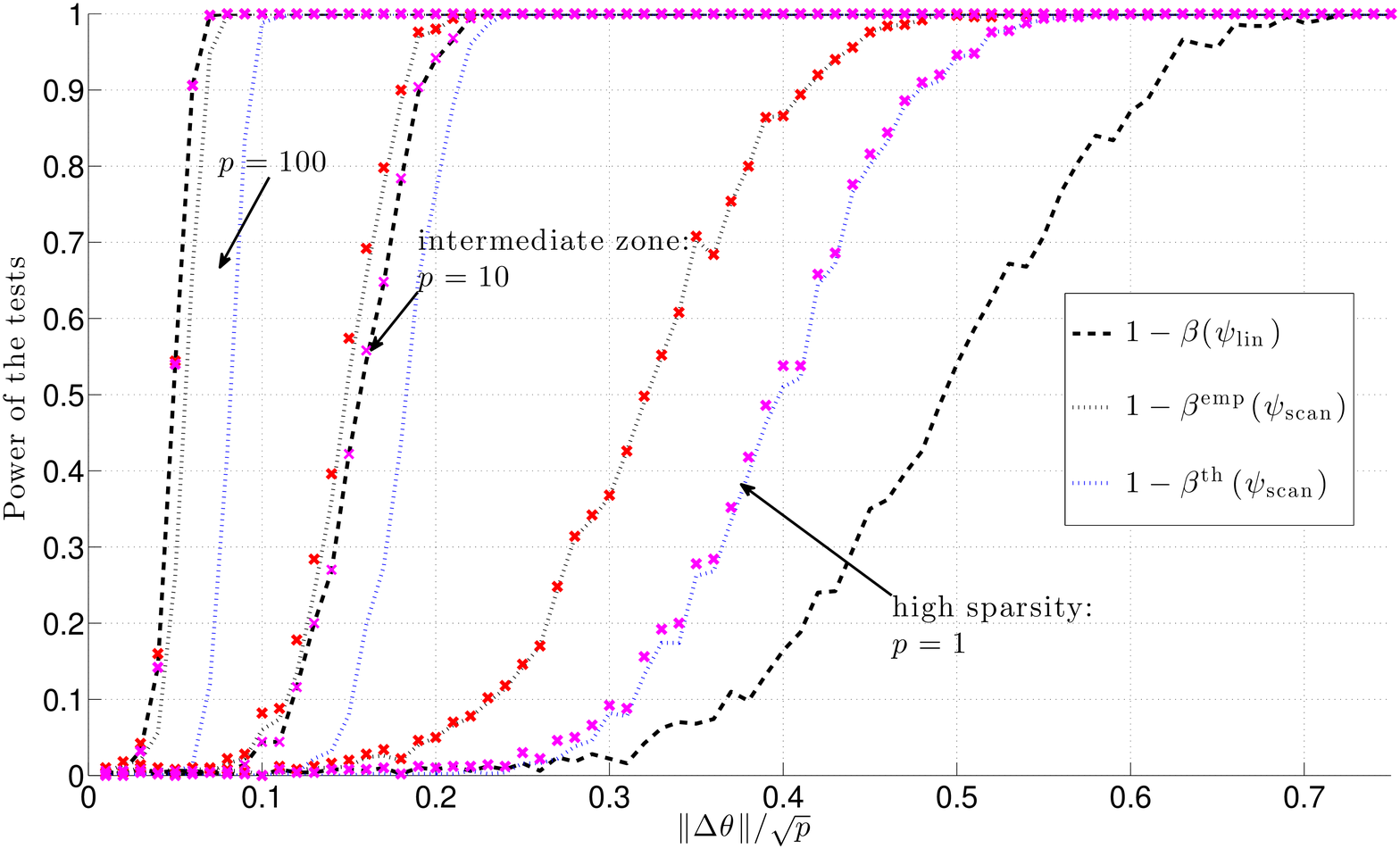}
\caption{Power of the test for $n=1000$, $d=100$, $\tau=n/2$ for $p=1,10,100$. In red and magenta colors the power of the proposed test $\psi^*=\psi_{\mathrm{lin}}\vee \psi_{\mathrm{scan}}$ is presented for the empirical and theoretical thresholds, respectively.}
\label{fig:3}
\end{figure}

\begin{figure}[htbp!]
 \centering
 \includegraphics[width=0.8\textwidth]{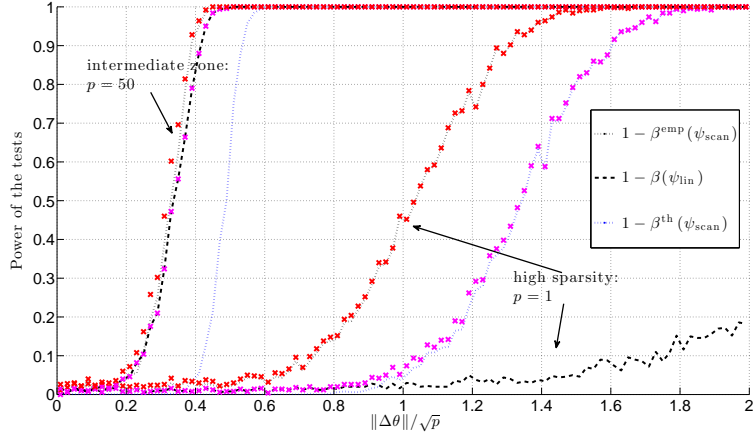}
 \caption{Power of the test for $n=100$, $d=1000$, $\tau=n/2$ for $p=1,50$. In red and magenta colors the power of the proposed test  $\psi^*$ is presented for the empirical and theoretical thresholds, respectively.}
 \label{fig:4}
 \end{figure}

 We also compared the power of the tests with different calibration strategies. The significance level was fixed to $0.05$. 
 A first calibration strategy consists in calibrating the scan test using the simulated quantiles instead of the threshold $T_p$. 
 A second calibration strategy consists in calibrating the test using the theoretical quantiles $T_p=\frac{\varkappa}{\sqrt {2p}}\Bigl[p\log(de/p)+\log(nd/\alpha)\Bigr]$.  However, the ``theoretical choice'' of the constant $\varkappa\approx 6.6$ obviously did not perform well in the simulations.  From the simulations we performed, we realized that choosing $\varkappa = 2$ gives much better results and seems to be robust to different values of $n$ and $d$.  In Figures~\ref{fig:3} and~\ref{fig:4}, we present the results of simulations for $n=1000$, $d=100$ and for $n=100$, $d=1000$.  In both cases, the change occurs in the middle of the interval of observations.  The linear statistic is calibrated using the quantile of level $\alpha/(nd)$ of $\chi^2$ distribution. We considered the cases of high sparsity, $p=1$, no sparsity, $p=d$ and the intermediate case $p\approx \sqrt d$. 
 
From Figure~\ref{fig:3} we see that the theoretical threshold $T_p$ gives quite good results. In the intermediate and medium sparsity cases the powers of the test (in blue solid and dashed lines for the empirical and the theoretical thresholds, respectively) are very close to each other. 
 
For the case of large $d$ on Figure~\ref{fig:4} we can see that the linear test is much less powerful in the situation of high sparsity than in the case of small $d$ as in Figure~\ref{fig:3}. 
  \begin{figure}[htbp!]
  \centering
  \includegraphics[width=0.8\textwidth]{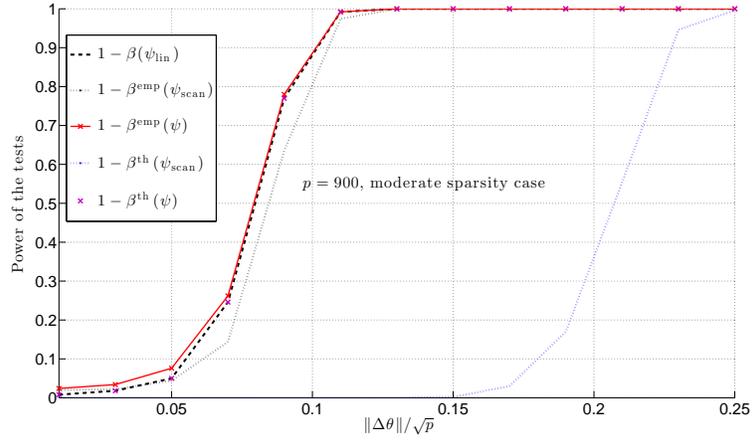}
  \caption{Power of the test for $n=100$, $d=1000$, $\tau=n/2$ for $p=900$. In red and magenta colors the power of the proposed test is presented for the empirical and theoretical thresholds, respectively.}
  \label{fig:5}
  \end{figure}
Figure~\ref{fig:5} shows the difference in the power of the tests for the case of moderate sparsity, $p=900$, $d=1000$ and $n=100$. The linear test outperforms the  scan test.

The last point we would like to mention is that the adaptive linear test can be calibrated by using of the empirical quantiles of the limiting process. We compared the accuracy of the theoretical quantile for the adaptive linear test given in~(\ref{thresholds2}). Note that this quantile depends on the step parameter $\eps\asymp m\sqrt{\log d/d}$. The results of the simulations show that choosing $m=2$ gives the best accuracy of the theoretical quantile of the adaptive linear test.

\bibliography{CP_sparse_test}

 \appendix
 \section{Upper bounds}
In what follows we use the following notation that takes origin from Lemma~\ref{Spok}, $\varkappa$ and $\varphi(t)=t-\log(1+t)$, $t\ge 0$

First, using~(\ref{Zj}) we reduce our model~(\ref{channels}) to the model of observations
\begin{equation}\label{Zj_obs}
Z_n(s)=\Delta\theta_\tau \mu_n (s)+\xi_n(s),\quad s=1,\dots,n-1,
\end{equation}
where 
\begin{equation}\label{mu}
\mu_n(s) =\sqrt{\frac{s(n-s)}n} \biggl(\frac{n-\tau }{n-s} \1\{s\le\tau\}
+ \frac{\tau}{s}\1\{s>\tau\} \biggr).
\end{equation}
and $\xi_n(s)=(\xi_n^1(s),\dots, \xi_n^d(s))^T$ are Gaussian vectors with  $\cN(0,1)$ dependent components given by
\begin{equation}\label{xi_jk}
\xi_n^j(s)=\sqrt{\frac{n}{s(n-s)}}\left(\sum_{i=1}^s \xi_i^j-\frac sn \sum_{i=1}^n \xi_i^j\right),\quad j=1,\dots,d.
\end{equation}
Note that the process $\xi_n^j(s)$, $s=1,\dots,n-1$ is a  discrete normalised Brownian bridge. 

\subsection{Upper bound for Problem (P1)}

The proof of Theorem~\ref{up_thm} is bases on two following lemmas.
\begin{lemma}\label{type_Ierr} 
The type I error in Problem~(P1) is bounded by $\alpha(\psi_p^*)\le \alpha(\psi_{\mathrm{lin}})+\alpha(\psi^p_{\mathrm{scan}})$, where
\begin{equation}\label{typeI_lin}
\alpha(\psi_{\mathrm{lin}})\le \begin{cases}
 e^{-\frac {\sqrt{2d}H}{\varkappa}},& H\ge \sqrt{d/2}\\
 e^{-\frac {2H^2}{\varkappa}},& H\le \sqrt{d/2};
 \end{cases}
 \end{equation}
 and
\begin{equation}\label{typeI_scan}
\alpha(\psi^p_{\mathrm{scan}})\le {d\choose p}  e^{-\frac{\sqrt{2p}T_p}\varkappa},\ T_p\ge\sqrt{p/2}.
\end{equation}
\end{lemma}
\begin{remark}\label{thrA}
 The thresholds 
\begin{equation}\label{thresholds1}
H=\begin{cases}
\frac{\varkappa}{\sqrt{2d}}\log\frac 1\alpha,& \alpha\le e^{-d/\varkappa}\\
\sqrt{\frac{\varkappa}2\log\frac 1\alpha},& \alpha\ge e^{-d/\varkappa}
\end{cases}
\quad \mbox{and}\quad
T_p =\frac\varkappa {\sqrt{2p}} \log \Bigl[{d\choose p}\frac 1\alpha\Bigr]
\end{equation}
define the tests $\psi_{\mathrm{lin}}$ and $\psi_{\mathrm{scan}}^p$ of significance level $\alpha$ so that $\alpha(\psi)\le 2\alpha$.
\end{remark}

{\bf Proof}.
Denote by $\chi^2_k$ a chi-squared random variable with $k$ degrees of freedom.
By Lemma~\ref{Spok},
\begin{align*}
\P_{\mathrm{H_0}}\{\psi_{\mathrm {lin}}=1\} &=\P_{\mathrm{H_0}}\left\{ \frac{\|Z_n(\tau)\|^2-d}{\sqrt{2d}}>H\right\}\\
&=\P\left\{ \frac{\|\xi_n(\tau)\|^2-d}{\sqrt{2d}}>H\right\}=\P \left\{ \chi^2_d> d+H\sqrt{2d} \right\}\\
&\le \exp\Bigl\{-\frac d2\varphi\Bigl(H\sqrt{\frac 2d}\Bigr)\Bigr\}
 \le \begin{cases}
 e^{-\frac {\sqrt{2d}H}{\varkappa}},& H\ge \sqrt{d/2}\\
 e^{-\frac {2H^2}{\varkappa}},& H\le \sqrt{d/2}.
 \end{cases}
\end{align*}
 The choice of $H$ as in (\ref{thresholds1}) provides $\P_{\mathrm{H_0}}\{\psi_{\mathrm {lin}}=1\}\le\alpha$. 
 
Let us now turn to the scan statistic. We have, again by Lemma~\ref{Spok},
\begin{align*}
\P_{\mathrm{H_0}} \{\psi_{\mathrm {scan}}=1\} &\le \P_{\mathrm{H_0}} \left\{\max_{m\in\cM(d,p)} \frac1{\sqrt{2p}} \Bigl(\|\Pi_m Z_n(s)\|^2-p\Bigr)>T_p\right\}\\
&= \P \left\{\max_{m\in\cM(d,p)} \frac1{\sqrt{2p}} \Bigl(\|\Pi_m \xi_n(s)\|^2-p\Bigr) > T_p\right\}\\
&\le  {d\choose p} \P \Bigl\{ \chi^2_p>p+T_p\sqrt{2p}\Bigr\}\\
&\le  {d\choose p} \exp\Bigl\{-\frac p2\varphi\Bigl(T_p\sqrt{\frac 2p} \Bigr)\Bigr\}\le {d\choose p}  e^{-\frac{\sqrt{2p}T_p}\varkappa},
\end{align*}
if $T_p\ge \sqrt{p/2}$. Let  $T_p =\frac\varkappa {\sqrt{2p}} \log \Bigl[{d\choose p}\frac 1\alpha\Bigr]$,  then $T_p\ge \sqrt{p/2} $ for any $\alpha>0$. For this value of the threshold $\P_{\mathrm{H_0}} \{\psi_{\mathrm {scan}}=1\} \le\alpha$.
The lemma is proved.
\hfill $\blacksquare$

\begin{lemma}\label{type_IIerr}
Assume that $p\in\{1,\dots,d\}$ and $\tau\in I$ are given, where $|I|\subseteq\{1,\dots,n-1\}$.
\begin{enumerate}
\item For any 
$\displaystyle{H\le \Bigl(\sqrt 2-1-c_1\Bigr) \frac{ r^2 nh(\tau)}{\sqrt {2d}}}$, where $c_1\in (0,\sqrt 2-1)$, 
the type II error for the linear test decreases exponentially as follows
\begin{equation}\label{typeII_lin}
\beta(\psi_{\mathrm {lin}},\Theta_p[r],\tau) \le \sqrt e\exp\left\{-c_1\frac{r^2 nh(\tau)}{\sqrt{2d}} \right\}.
\end{equation}
\item Let  $\displaystyle{c_2\equiv c_2(p,\tau,r)=(1+\sqrt{2/p})^{-1}\frac{r^2 nh(\tau)}{T_p\sqrt{2p}}}$.
Then the type II error for the scan test is bounded by
\begin{equation}\label{typeII_scan}
\beta (\psi_{\mathrm {scan}}^p,\Theta_p[r],\tau) \le \sqrt e\exp\left\{-(c_2-1) T_p\right\}.
\end{equation}
\end{enumerate}
\end{lemma}

{\bf Proof}. 
Let $m$ be the true subset of $\cM(d,p)$ of indexes of components of the mean with a change. For the scan statistic we have
\begin{align*}
\P_{\Delta\theta_\tau}\Bigl\{L^p_{\mathrm {scan}}(\tau)<T_p\Bigr\}&\le \min_{m\in\cM(d,p)}\P_{\Delta\theta_\tau} \left\{ 
\frac1{\sqrt{2p}}\biggl(\|\Pi_m Z_n(\tau)\|^2-p\biggr)<T_p\right\}\\
&\le \P_{\Delta\theta_\tau}\biggl\{ \frac1{\sqrt{2p}}\biggl(\|\Pi_m Z_n(\tau)\|^2-p\biggr)<T_p\biggr\}\\
&\le e^{T_p+\sqrt{p/2}} \E_{\Delta\theta_\tau} \exp \Bigl(-\frac1{\sqrt{2p}}\|\Pi_m Z_n(\tau)\|^2\Bigr).
\end{align*}
Denote $\lambda=-1/\sqrt{2p}$. We have for $h(\tau)$ defined in~(\ref{htau})
\begin{align*}
\E_{\Delta\theta_\tau}  \exp \Bigl(\lambda\|\Pi_m Z_n(\tau)\|^2\Bigr) &= \exp\biggl(\frac{\lambda \|\Delta\theta_\tau\|^2nh(\tau)}{1-2\lambda}-\frac{p}2\log(1-2\lambda)\biggr)\\
&\le \exp\biggl(-\frac{\|\Delta\theta_\tau\|^2nh(\tau)}{\sqrt{2p}}(1+\sqrt{2/p})^{-1} -\frac{p}2(\sqrt{2/p}-1/p)\biggr)\\
&= \exp\biggl(-\frac{\|\Delta\theta_\tau\|^2nh(\tau)}{\sqrt{2p}}(1+\sqrt{2/p})^{-1} -\sqrt{p/2}+1/2\biggr).
\end{align*}
Therefore,
$$
\P_{\Delta\theta_\tau}\Bigl\{L^p_{\mathrm {scan}}(\tau)<T_p\Bigr\}\le \exp\biggl(T_p  -\frac{\|\Delta\theta_\tau\|^2nh(\tau)}{\sqrt{2p}}(1+\sqrt{2/p})^{-1} +1/2\biggr)
$$
which implies
\begin{align*}
\beta(\psi_{\mathrm {scan}}^p,\Theta_p[r],\tau)&=\sup_{\Delta\theta_\tau\in\Theta_p[r]} \P_{\Delta\theta_\tau}\{L^p_{\mathrm {scan}}(\tau)<T_p\}\\
&\le  \exp\biggl(T_p  -\frac{r^2nh(\tau)}{\sqrt{2p}}(1+\sqrt{2/p})^{-1} +1/2\biggr)\\
&=\sqrt e\exp\left\{-(c_2-1) T_p\right\}.
\end{align*}

For the linear statistic we can proceed exactly as in the case of the scan statistic. We have
$$
\P_{\Delta\theta_\tau}\{L_{\mathrm {lin}}(\tau)<H\}\le \exp\biggl(H -\frac{\|\Delta\theta_\tau\|^2nh(\tau)}{\sqrt{2d}}(1+\sqrt{2/d})^{-1} +1/2\biggr).
$$
Therefore, for any $c_1\in(0,\sqrt 2-1)$ such that $(1+\sqrt{2/d})^{-1}- \frac{H\sqrt{2d} }{\|\Delta\theta_\tau\|^2nh(\tau)}> c_1$ 
$$
\beta(\psi_{\mathrm {lin}},\Theta_p[r],\tau)\le \sqrt e\exp\biggl\{-\Bigl( (1+\sqrt{2/d})^{-1}-c_1\Bigr)\frac{r^2 nh(\tau)}{\sqrt{2d}}\biggr),
$$
and the lemma follows. \hfill$\blacksquare$

{\bf Proof of Theorem~\ref{up_thm}}. 
Lemma~\ref{type_Ierr}  implies $\alpha(\psi_p^*)\to 0$ as $d\to\infty$ for any $H>\sqrt{(\varkappa/2)\log d}\to\infty$ and for $T_p$  defined in (\ref{Tp}). 
On the other hand, 
$$
\beta (\psi_p^*,\Theta_p[r],\tau)\le \min\Bigl(\beta(\psi_{\mathrm {lin}},\Theta_p[r],\tau),\beta(\psi^p_{\mathrm {scan}},\Theta_p[r],\tau)\Bigr).
$$
If~(\ref{H}) is satisfied, then there exists $c_1\in(0,\sqrt 2-1)$ such that the conditions of the first part of Lemma~\ref{type_IIerr} hold.
Therefore, by Lemma~\ref{type_IIerr}
$$
\beta(\psi_{\mathrm {lin}},\Theta_p[r],\tau)\le \sqrt e\exp\left\{-c_1\frac{r^2 nh(\tau)}{\sqrt{2d}} \right\}.
$$
Since (\ref{up_adapt_lin}) holds true and $c_1>0$, we have $\beta(\psi_{\mathrm {lin}},\Theta_p[r],\tau)\to 0$ as $d\to\infty$. 

By Lemma~\ref{type_IIerr}
$$
\beta (\psi_{\mathrm {scan}},\Theta_p[r],\tau) \le \sqrt e\exp\left\{-(c_2-1) T_p\right\}.
$$
where $T_p=\frac \varkappa {\sqrt{2 p}}\log \biggl[(1+\delta) {d\choose p}\biggr]$ as defined in~(\ref{Tp}). 
Recall the following fact: $p\log (d/p)\le \log {d\choose p}\le p\log (de/p)$. We have
$$
\sqrt{2p}T_p=\varkappa \Bigl(\log{d\choose p} +\log (1+\delta)\Bigr) <\varkappa p\biggl(\log\frac{de}p +\frac{\log (1+\delta)}p\biggr).
$$
Then we have
\begin{align*}
c_2(p,\tau,r) &= (1+\sqrt{2/p})^{-1}\frac{r^2 nh(\tau)}{\sqrt{2p}T_p}\\
&>  \frac{r^2nh(\tau)}{\varkappa p\log (d/p)}\cdot\frac{(1+\sqrt{2/p})^{-1}}{1+[\log(d/p)]^{-1}+[p\log (d/p)]^{-1}\log(1+\delta)}.
\end{align*}
Since $p\to\infty$ and $\log(d/p)\to\infty$  as $d\to \infty$, by condition~(\ref{scan_dbp}) of Theorem~\ref{up_thm} we have
$$
\liminf_{d\to\infty}c_2(p,\tau,r)=\liminf_{d\to\infty} \frac{(1+\sqrt{2/p})^{-1}r^2nh(\tau)}{\sqrt{2p}T_p}>\liminf_{d\to\infty} \frac{r^2nh(\tau)}{\varkappa p\log (d/p) }>1.
$$
Let $\liminf\limits_{d\to\infty} c_2(p,\tau,r)\ge c_0>1$. This implies
\begin{align*}
\beta(\psi^p_{\mathrm {scan}},\Theta_p[r],\tau)&\le \sqrt e\exp\Bigl\{-(c_0-1) T_p\Bigr\}\\
&\le  \sqrt e\exp\biggl\{-\varkappa(c_0-1) \Bigl(\sqrt{\frac p2} \log\frac dp +\frac{\log(1+\delta)}{\sqrt{2p}}\Bigr) \biggr\}
\to 0,\quad d\to\infty.
\end{align*}
This gives $\gamma(\psi_p^*,\Theta_p[r],\tau)\to 0$ as $d\to\infty$.\hfill $\blacksquare$

\subsection{Upper bound for Problem (P2)}
The proof of the upper bound in the adaptive case is similar to the proof in the case of known $p$ and $\theta$. 
It is also based on two lemmas. 
\begin{lemma}\label{type_Ierr_adapt} 
Let $\eps>0$, $\eps/n\to 0$ as $n\to\infty$. Then $\alpha(\psi^*)\le \alpha(\psi_{\mathrm {lin}}) +\alpha(\psi_{\mathrm {scan}})$, where
$$
\alpha(\psi_{\mathrm{lin}}) \le  2 \log_{1+\eps} n  \begin{cases}
e^{-\frac 1{\varkappa(1+\eps)}\sqrt{2d}(H-\eps\sqrt{d/2})},& H\ge (1+2\eps)\sqrt{d/2}\\
 e^{-\frac 2{\varkappa(1+\eps)^2}(H-\eps\sqrt{d/2})^2} ,&H\le (1+2\eps)\sqrt{d/2}.
\end{cases}
$$
and
$$
\alpha(\psi_{\mathrm{scan}})\le  n\sum_{p=1}^d  {d\choose p}  e^{-\frac{\sqrt{2p}T_p}\varkappa}.
$$ 
\end{lemma}
\begin{remark}\label{rem:approx_thresh}
The thresholds
\begin{equation}\label{thresholds2}
H=\begin{cases}
\frac{\varkappa(1+\eps)}{\sqrt{2d}}\log \Bigl[\frac{\log_{1+\eps}n}\alpha \Bigr]+\eps\sqrt{\frac d2},& \alpha\le \log_{1+\eps}n e^{-d/\varkappa}\\
\sqrt{\frac{\varkappa(1+\eps)^2}{2}\log\Bigl[\frac{\log_{1+\eps}n}\alpha \Bigr]}+\eps\sqrt{\frac d2},& \alpha\ge \log_{1+\eps}n e^{-d/\varkappa}
\end{cases}
\end{equation}
and $T_p =\frac\varkappa {\sqrt{2p}} \log \Bigl[{d\choose p}\frac {nd}\alpha\Bigr]$ define two tests of significance level $\alpha$ each such that $\alpha_B(\psi^*)\le 2\alpha$.
\end{remark}

{\bf Proof}. Let us start with the linear test. We have 
\begin{align*}
\P_{\mathrm{H_0}}\{\psi_{\mathrm {lin}}=1\} &=\P_{\mathrm{H_0}}\Bigl\{ \max_{s=1,\dots,n-1}\frac{\|Z_n(s)\|^2-s}{\sqrt{2d}}>H\Bigr\}\\
 &=\P \Bigl\{ \max_{s=1,\dots,n-1} \|\xi_n(s)\|^2 >d+ H\sqrt{2d}\Bigr\}
\end{align*}
Note that $\xi_n^j(s)\deq B_j(t)$ for $t=s/n$, where $B_j(t)$ are standard independent Brownian brigdes on $[0,1]$. To bound the probability on the last display we shall first pass to the Brownian bridge $B_j$ and 
then to the standard Wiener process $W_j$ using the time change $u=t/(1-t)$,
\begin{align*}
\P \biggl\{ \max_{s=1,\dots,n-1}\|\xi_n(s)\|^2>d+ H\sqrt{2d}\biggr\}
&\le \P \biggl\{ \sup_{\frac 1n\le t\le 1-\frac1n }\sum\limits_{j=1} ^d \frac{B_j^2(t)}{t(1-t)}>d+ H\sqrt{2d}\biggr\}\\
&=  \P \biggl\{ \sup_{\frac 1{n-1}\le u\le n-1}\sum\limits_{j=1} ^d \frac{W_j^2(u)}u>d+ H\sqrt{2d}\biggr\}.
\end{align*}
Proceeding like in~\cite{Vostrikova:1981} we split the interval $[1,n-1]$ into $N(\eps)=\log_{1+\eps} (n-1)$ blocks $[a_k,a_{k+1}]$, where $a_k=(1+\eps)^k$. Then we estimate the probability of maximum over each block. Denote $x^2=d+H\sqrt{2d}$. We have
\begin{align*}
\P  &\biggl\{\sup_{u\in [\frac1{n-1}, n-1]}\sum\limits_{j=1} ^d \frac{W_j^2(u)}u>x^2\biggr\}\\
 &\le \P\biggl\{\sup_{u\in [\frac1{n-1},1]} \frac{\|W(u)\|}{\sqrt u}> |x| \biggr\} 
 +\sum_{k=0}^{N(\eps)-1} \P\biggl\{\sup_{s\in [a_k,a_{k+1}-1]}  \frac{\|W(u)\|}{\sqrt u}> |x| \biggr\}  \\
 &\le \P\biggl\{\sup_{u\in[0,1]} \|W(u)\|>\sqrt{n-1}|x|\biggr\}\\
  &\qquad\qquad+\sum_{k=0}^{N(\eps)-1} \P\biggl\{\sup_{u\in [a_k,a_{k+1}-1]} \|W(u)\|> (1+\eps)^{k/2}|x| \biggr\}.
\end{align*}
Next, using the scaling property of the Wiener process we have
$$
\P\Bigl\{\sup_{u\in [a_k,a_{k+1}-1]} \|W(u)\|> (1+\eps)^{k/2} x\Bigr\} = \P\Bigl\{\sup_{u\in [1,1+\eps-a_k^{-1}]} \|W(u)\|> x\Bigr\}.
$$
Therefore,
\begin{align*}
\P&\biggl\{\sup_{u\in [\frac1{n-1}, n-1]} \frac{\|W(u)\|}{\sqrt u}> |x|\biggr\}\\
&\le \P\Bigl\{\sup_{u\in [0,1]} \|W(u)\|> \sqrt{n-1}|x|\Bigr\}+ (N(\eps)-1) \P\Bigl\{\sup_{u\in [1,1+\eps]} \|W(u)\|> x\Bigr\}\\
&\le \log_{1+\eps} n \P\Bigl\{\sup_{u\in [0,1+\eps]} \|W(u)\|> x\Bigr\} \le 2 \log_{1+\eps} n \P\Bigl\{\|\xi\|^2>\frac{x^2}{1+\eps}\Bigr\},
\end{align*}
where $\xi$ is the standard Gaussian $d$-dimensional vector. 
Using the bound from Lemma~\ref{Spok} we have
\begin{align*}
\P_{\mathrm{H_0}}\{\psi_{\mathrm{lin}}=1\}&\le \P \biggl\{ \sup_{u\in [\frac1{n-1}, n-1]}\sum\limits_{j=1} ^d \frac{W_j^2(u)}u>d+\sqrt{2d}H\biggr\}\\
&\le  2 \log_{1+\eps} n  \exp\biggl\{-\frac d2\varphi\biggl(\frac{H}{1+\eps}\sqrt{\frac 2d}-\frac \eps{1+\eps}\biggr)\biggr\}\\
&\le 2 \log_{1+\eps} n  \begin{cases}
e^{-\frac 1{\varkappa(1+\eps)}\sqrt{2d}(H-\eps\sqrt{d/2})},& H\ge (1+2\eps)\sqrt{d/2}\\
 e^{-\frac 2{\varkappa(1+\eps)^2}(H-\eps\sqrt{d/2})^2} ,&H\le (1+2\eps)\sqrt{d/2}.
 \end{cases}
 \end{align*}
For $H$ chosen as in~(\ref{thresholds2}) we have $\P_{\mathrm{H_0}}\{\psi_{\mathrm{lin}}=1\}\le\alpha$.

Let us now turn to the scan statistic. Here we can just use the union bound to show the correct rates. We have like in the case of Problem~(P1)
\begin{align*}
\P_{\mathrm{H_0}}\{\psi_{\mathrm {scan}}=1\} &= \P_{\mathrm{H_0}}\left\{ \max_{s=1,\dots,n-1} L_{\mathrm {scan}}(s)>1\right\}\\
 &\le n \sum_{p=1}^ d\P \left\{\max_{m\in\cM(d,p)} \frac1{\sqrt{2p}} \left(\|\Pi_m \xi_n(s)\|^2-p\right) > T_p\right\}\\
 &\le  n\sum_{p=1}^d {d\choose p} \P \left\{ V_p>p+T_p\sqrt{2p}\right\}\\
 &\le  n\sum_{p=1}^d {d\choose p} \exp\left\{-\frac p2\varphi\left(T_p\sqrt{\frac 2p} \right)\right\}
 \le n\sum_{p=1}^d  {d\choose p}  e^{-\frac{\sqrt{2p}T_p}\varkappa}.
\end{align*}
In particular,  $T_p$ chosen as specified in Remark~\ref{rem:approx_thresh} provides us with the test of significance level $\alpha$. \hfill $\blacksquare$

\begin{lemma}\label{type_IIerr_adapt} 
Let $I\subseteq\{1,\dots,n-1\}$. Then 
\begin{enumerate}
\item For any 
$\displaystyle{H\le \left(\sqrt 2-1-c_1\right) \min_{\tau\in I} \frac{ r^2 nh(\tau)}{\sqrt {2d}}}$, where $c_1\in (0,\sqrt 2-1)$, 
the type II error for the linear test decreases exponentially as
\begin{equation}\label{typeII_lin_adapt}
\beta^* (\psi_{\mathrm {lin}},\Theta[r],I) \le \sqrt e\exp\left\{-c_1\min_{\tau\in I}\frac{r^2 nh(\tau)}{\sqrt{2d}} \right\}.
\end{equation}
\item Let  $\displaystyle{c_2(p,\tau,r)=(1+\sqrt{2/p})^{-1}\frac{r^2 nh(\tau)}{T_p\sqrt{2p}}}$.
Then the type II error for the scan test is bounded by
\begin{equation}\label{typeII_scan_adapt}
\beta^* (\psi_{\mathrm {scan}},\Theta[r],I)\le \sqrt e\exp\left\{-\min_{p=1,\dots,d}\min_{\tau \in I}\left[(c_2(p,\tau,r)-1) T_p\right]\right\}.
\end{equation}
\end{enumerate}

\end{lemma}
{\bf Proof}. 
For the linear test we have
\begin{align*}
\beta^*(\psi_{\mathrm{lin}},\Theta[r],I) &=\sup_{p=1,\dots,d}\sup_{\tau\in I}\sup_{\Delta\theta_\tau\in\Theta_p[r]} \P_{\Delta\theta_\tau}\{\psi_{\mathrm{lin}}=1\}\\
&=\sup_{p=1,\dots,d}\sup_{\tau\in I} \beta(\psi_{\mathrm{lin}},\Theta_p[r],\tau)\\
&\le  \sup_{\tau\in I}  \sqrt e \exp\left\{-c_1\frac{r^2nh(\tau)}{\sqrt{2d}}\right\}.
\end{align*}
where the last line follows from Lemma~\ref{type_IIerr}. 
The result for the scan test also follows from a direct application of Lemma~\ref{type_IIerr}. We have
\begin{align*}
\beta^*(\psi_{\mathrm{scan}},\Theta[r],I) &=\sup_{p=1,\dots,d}\sup_{\tau\in I}\sup_{\Delta\theta_\tau\in\Theta_p[r]} \P_{\Delta\theta_\tau}\{\psi_{\mathrm{scan}}=1\}\\
&=\sup_{\tau\in I} \sup_{p=1,\dots,d}\beta(\psi_{\mathrm{scan}},\Theta_p[r],I)\\
&\le  \sup_{\tau\in I}  \sup_{p=1,\dots,d}\exp\left(T_p  -\frac{r^2nh(\tau)}{\sqrt{2p}}(1+\sqrt{2/p})^{-1} +1/2\right)\\
&\le \sqrt e\exp\left\{-\min_{p=1,\dots,d}\min_{\tau \in I}\left[(c_2(p,\tau,r)-1) T_p\right]\right\}.
\end{align*}
The lemma is proved. \hfill $\blacksquare$

{\bf Proof of Theorem~\ref{up_thm_adapt}}.
 Let $\eps=m\sqrt{\log d/d}$. Note that for any $m\ge 1$ and $d\ge 3$
 \begin{align*}
 -\log\log&\left(1+m\sqrt{\frac{\log d}d}\right)\le -\log\left(m\sqrt{\frac{\log d}d}\left(1-\frac m2\sqrt{\frac{\log d}d}\right)\right)\\
 &=-\log (m\sqrt{\log d})-\log(\sqrt d-(m/2)\sqrt{\log d})+\log d\le \log d-\log m. 
 \end{align*}
 
Suppose that $H$ satisfies~(\ref{H_adapt}) and $c_0\sqrt{\log(d\log n)} \le H \le  (1+2\eps)\sqrt{d/2}$ for $c_0>\sqrt{\varkappa}$. The threshold $H$ can be always chosen in this way since~(\ref{up_adapt_lin}) holds.Then we have from Lemma~\ref{type_Ierr_adapt}
 \begin{align*}
 \alpha(\psi_{\mathrm{lin}}) &\le 2\log n\exp\left\{-\log\log(1+m\sqrt{\log d/d}) -\frac{2\left(c_0\sqrt{\log(d\log n)}-m\sqrt{\log d/2}\right)^2 }{\varkappa(1+m\sqrt{\log d/d})^2}\right\}\\
 &\le \frac 2m\exp\left\{\log\log n+\log d -\frac 2\varkappa \frac {\left((c_0-m)\sqrt{\log d/2}+c_0\sqrt{\log\log n/2}\right)^2 }{2(1+m)^2} \right\},
 \end{align*}
where the inequality $\sqrt{a+b}\ge \sqrt{a/2}+\sqrt{b/2}$ for $a,b>0$ was used in the last line.
Next, since $(a+b)^2>a^2+b^2$ for $a,b>0$, we obtain
 \begin{align*}
 \alpha(\psi_{\mathrm{lin}})  &\le \frac 2m  \exp\left\{\log\log n+\log d  -\frac {(c_0-m)^2\log d+c_0^2\log\log n}{\varkappa (1+m)^2} \right\}\\
&\le \frac 2m\exp\left\{\left(1-\frac{c_0^2}{\varkappa(1+m)^2}\right)\log\log n +\left(1-\frac{(c_0-m)^2}{\varkappa(1+m)^2}\right)\log d \right\}.
 \end{align*}
Choose $m\in (0,(c_0-\sqrt{\varkappa})/(\sqrt{\varkappa}+1)$, where $c_0>\sqrt{\varkappa}$. For such a  choice of $m$ we have
 $(c_0-m)^2>\varkappa(1+m)^2$ and, consequently, $ \alpha_B(\psi_{\mathrm{lin}})\to 0$ as~$d\to\infty$.  
 
Since $H$ satisfies (\ref{H_adapt}), from Lemma~\ref{type_IIerr_adapt} it follows immediately that\linebreak $\beta^*(\psi_{\mathrm{lin}},\Theta[r],I)\to 0$. 
 
Note also that from Lemma~\ref{type_Ierr_adapt} it follows that for any  $H\ge (1+2\eps)\sqrt{d/2}$ 
\begin{align*}
\alpha(\psi_{\mathrm{lin}}) &\le 2\log n\exp\left\{-\log\log(1+m\sqrt{\log d/d}) -\frac{\sqrt{2d}(1+m\sqrt{\log d/2})}{\varkappa(1+m\sqrt{\log d/d})}\right\}\\
&\le \frac{2}m\exp\left\{\log\log n+\log d -\frac {\sqrt d}{\varkappa} \left(1+m\sqrt{\log d/2}\right)\right\}
\to 0,\quad d\to\infty
\end{align*}
if $\log n/d \to 0 $ as $d\to\infty$. The latter asymptotics always holds if (A\ref{asymp_cond2}) is satisfied since $d/p\ge \log(d/p)$ for all $d\ge 2$. Note that in this case
$$
H\ge(1+2\eps)\sqrt{d/2}=\sqrt{d/2}+m\sqrt{2\log d}>\sqrt{d/2}.
$$
This choice of $H$ provides a faster rate in the type~I error if the change in the signal is sufficiently large, $\|\Delta\theta_\tau\|^2\asymp d$.

Let us turn to the scan test. For $T_p=\frac{\varkappa}{\sqrt{2p}} \Bigl[\log{d\choose p}+\log(ndp^2)\Bigr]$, where $\delta>0$ we have $T_p\ge \sqrt{p/2}$. Thus, applying Lemma~\ref{type_Ierr_adapt} we obtain
\begin{align*}
\alpha(\psi_{\mathrm{scan}})&\le n\sum_{p=1}^d  {d\choose p}  e^{-\frac{\sqrt{2p}T_p}\varkappa}\\
&= n\sum_{p=1}^d  {d\choose p}  \exp\Bigl\{-\varkappa\Bigl[\log{d\choose p}-\log(ndp^2)\Bigr]\Bigr\}\\
&=\frac 1d\sum_{p=1}^d \frac 1{p^2} \le \frac 1d \frac{\pi^2}6 \to 0,\quad d\to\infty.
\end{align*}

For the type II error we have from Lemma~\ref{type_IIerr_adapt}
$$
\beta^*(\psi_{\mathrm{scan}},\Theta[r],I) 
\le \sqrt{e} \exp\Bigl\{-\min_{p=1,\dots,d}\min_{\tau\in I} [(c_2(p,\tau,r)-1)T_p]\Bigr\},
$$
where
\begin{align*}
c_2(d,p,r) &= (1+\sqrt{2/p})^{-1}\frac{r^2 nh(\tau)}{\sqrt{2p}T_p}\\
&= (1+\sqrt{2/p})^{-1}\frac{r^2 nh(\tau)}{\varkappa \log{d\choose p} +\varkappa \log(ndp^2)}\\
&> (1+\sqrt{2/p})^{-1}\frac{r^2 nh(\tau)}{\varkappa p\log(de/p)+\varkappa \log(nd)+2\varkappa \log p}\\
&=  \frac{r^2nh(\tau)}{\varkappa p\log (d/p)}\cdot\frac{(1+\sqrt{2/p})^{-1}}{1+[\log(d/p)]^{-1}+[\log(nd)+2\log p]\cdot[p\log(d/p)]^{-1}}.
\end{align*}
By Assumptions~(A\ref{asymp_cond1}--A\ref{asymp_cond2}), we have $\log(d/p)\to\infty$  , $p\to\infty$, and $\log (nd)/(p\log d)\to 0$ as $d\to \infty$. Therefore, by condition~\ref{up_adapt_scan} of Theorem~\ref{up_thm_adapt} we have for all $\tau\in I$ and $p=1,\dots,d$
$$
\liminf_{d\to\infty} \min_{\tau\in I}\min_{p=1,\dots,d} c_2(d,p,r)>\liminf_{d\to\infty} \min_{\tau\in I}\min_{p=1,\dots,d} \frac{r^2nh(\tau)}{\varkappa p\log (d/p) }>1.
$$
Combining this result with the fact that $T_p\to \infty$ as $d\to\infty$ for all $p=1,\dots,d$ we obtain that $\beta^*(\psi_{\mathrm {scan}},\Theta[r],I)\to 0$ as  $d\to\infty$.
This proves the theorem.\hfill $\blacksquare$ 

\section{Proof of the lower bound}
We have to show that under conditions~of Theorem~\ref{low_thm}
$$
\inf_{\psi\in[0,1]} \gamma(\psi,\Theta_p[r],\tau)\to 1,\quad d\to\infty.
$$
The classical approach to the construction of lower bounds in testing problems goes back to the works of~\cite{Ingster&Suslina:2003}. 
We also refer to~\citep{Baraud:2002} for details on the non-asymptotic aspect of the problem.

Recall the new model~(\ref{Zj_obs}),  Gaussian process~(\ref{xi_jk}) and the function $\mu_n$ defined in~(\ref{mu}).
Denote by $ Z=(Z_n^j(s))$ the $d\times n-1$ matrix of observations~(\ref{Zj_obs}). 

Let $\widetilde \Theta_p[r]\subset \Theta_p [r]$ be a subset of $d$-dimensional vectors $v$ with the norm $\|v\|=r$:
$$
\widetilde \Theta_p[r]=\{v\in V_p^d:\ \|v\|=r\}.
$$
Let $\P_{\pi,r}$ be a prior distribution on $v\in\widetilde \Theta_p[r]$ that will be defined later. We have 
\begin{align*}
\inf_{\psi\in[0,1]}\gamma (\psi,\Theta_p[r],\tau) &=\inf_{\psi\in[0,1]} \Bigl(\alpha(\psi)+\beta(\psi,\Theta_p[r],\tau)\Bigr)\\
&\ge \inf_{\psi\in[0,1]} \Bigl(\alpha(\psi)++\beta(\psi,\widetilde\Theta_p[r],\tau)\Bigr) \\
&\ge 1-\frac12 \| \P_{\pi,r}-\P_0\|_{TV}\\
&\ge 1-\frac12\Bigl(\E_0 L^2_{\pi,r} (Z)-1\bigr)^{1/2},
\end{align*}
where $\|\cdot\|_{TV}$ denotes the total variation norm and 
$$
 L_{\pi,r} (Z)=\frac{d\P_{\pi,r}}{d\P_0} (Z)
$$
be the corresponding likelihood ratio. Thus to prove the lower bound  it is sufficient to show that
\begin{equation}\label{L2conv}
\limsup_{d\to\infty} \E_0[L^2_{\pi,r} (Z)]\to 1.
\end{equation}

In what follows in this section we assume that $m\in\cM(d,p)$ is uniformly distributed over $\cM(d,p)$ according to the measure $\pi$. 
Let $\P_{m,r}$ be the prior distribution on $\Delta\theta_\tau\in \widetilde \Theta_p[r]$, where $\Delta\theta_\tau$ is the vector of change in mean in components indexed by the set $m$. This prior will be specified later. Thus, we have the following mixture of priors on $\Delta\theta_\tau$ and $m$,
$$
\P_{\pi,r}={{d\choose p}}^{-1}\sum_{m\in\cM(d,p)} \P_{m,r}
$$
and the corresponding  likelihood ratio given by
$$
L_{\pi,r}(Z)={{d\choose p}}^{-1} \sum_{m\in\cM(d,p)} \frac{d\P_{m,r}}{d\P_0} (Z).
$$

\begin{lemma}\label{lemLRHG}
Assume that the prior distribution $\P_{m,r}$ on $\Delta\theta\in\widetilde \Theta_p[r]\subset V_p^d$ is given by 
$$
\Delta\theta_j=\frac{ r}{\sqrt p}\eps_j^m,\quad j=1,\dots,d.
$$
Let $Y$ be a binomial r.v. with parameters $p\in\{1,\dots,d\}$ and $p/d$. 
\begin{enumerate}[(i)]
\item If $\eps_j^m=\begin{cases}\pm1,& j\in m\\
 										0,& j\notin m
 										\end{cases}$,
then $\E_0 [L^2_{\pi,r} (Z)] \le \E  \Bigl[\cosh \Bigl(\frac{r^2}p nh(\tau)\Bigr)\Bigr]^{Y}$.
\item If  $\eps_j^m=\begin{cases}1,& j\in m\\
 										0,& j\notin m
 										\end{cases}$,
then $\E_0 [L^2_{\pi,r} (Z)] \le \E  \Bigl[\exp \Bigl(\frac{r^2}p nh(\tau)\Bigr)\Bigr]^{Y}$.
\end{enumerate}
\end{lemma}
{\bf Proof} follows a general structure originally proposed in~\cite{Baraud:2002}. 
Note that $\xi_n^j(s)$, $s=1,\dots,n-1$ in~(\ref{xi_jk}) are not independent. Denote
$\tilde\xi^j=(\xi_n^j(1),\dots \xi_n^j(n-1) )^T$. Let $A$ be a $(n-1)\times n$ matrix with components
$$
a_{st}=\sqrt{\frac{s(n-s)}{n}}\Bigl[\frac1s \1\{t\le s\}-\frac1{n-s}\1\{t>s\}\Bigr].
$$
Then $\tilde\xi^j=A\xi^j$, where $\xi^j=(\xi_{1j},\dots,\xi_{nj})^T\sim\cN(0,I_d)$, $j=1,\dots,d$.
Denote\linebreak $\mu_n=(\mu_n(1),\dots,\mu_n(n-1))^T$. Thus for each $m$ the likelihood ratio can be written in terms of matrix $A$ and function $\mu_n$,
\begin{align*}
\frac{d\P_{m,r}}{d\P_0} (Z)&= \prod_{j\in m} \frac{\exp\Bigl\{-\frac12(Z_j-\Delta\theta_j\mu_n)^T(AA^T)^{-1}(Z_j-\Delta\theta_j\mu_n)\Bigr\}}
{\exp\Bigl\{-\frac12(Z_j^T(AA^T)^{-1}Z_j)\Bigr\}}\\
&=\prod_{j\in m} \exp\Bigl\{ \Delta\theta_j \mu_n^T(AA^T)^{-1} Z_j-\frac12 \Delta\theta_j^2 \mu_n^T (AA^T)^{-1}\mu_n\Bigr\}\\
&= \exp\Bigl\{-\frac 12\|\Delta\theta\|^2\mu_n^T (AA^T)^{-1}\mu_n\Bigr\} \prod_{j\in m} \exp\Bigl\{ \Delta\theta_j \mu_n^T(AA^T)^{-1} Z_j\Bigr\}
\end{align*}
It can be shown  that the matrix $B=(AA^T) ^{-1}$ is tridiagonal with components
$$
b_{ij}=2i\Bigl(1-\frac in\Bigr)\1\{i=j\}-\frac1n \sqrt{ij(n-i)(n-j)}\1\{|i-j|=1\},\quad i,j=1,\dots,n-1.
$$
Direct calculation gives the following relation,
$$
\mu_n^T (AA^T)^{-1}\mu_n=\frac{\tau(n-\tau)}n.
$$
Denote 
$$
\rho^2= \|\Delta\theta\|^2nh(\tau).
$$
Then 
$$
\frac{d\P_{m,r}}{d\P_0} (Z) = e^{-\rho^2/2} \exp\biggl\{-\sum_{j\in m} \Delta\theta_j \mu_n^T(AA^T)^{-1} Z_j\biggr\}.
$$

First, let us consider the ``Rademacher prior'', that is the prior defined in~(i) of Lemma~\ref{lemLRHG}. We have
\begin{align*}
L^2_{\pi,\Delta\theta} (Z)&={{d\choose p}}^{-1} e^{-\rho^2/2} \sum_{m\in\cM(d,p)} \E_{\eps}   \exp\biggl\{-\sum_{j\in m} \frac{\|\Delta\theta\|}{\sqrt p}\eps_j \mu_n^T(AA^T)^{-1} Z_j\biggr\}\\
&= {{d\choose p}}^{-1} e^{-\rho^2/2} \sum_{m\in\cM(d,p)} \prod_{j\in m}\cosh\Bigl(\frac{\|\Delta\theta\|}{\sqrt p}\mu_n^T(AA^T)^{-1} Z_j\Bigr)
\end{align*}
Now we have to calculate $\E_0 [L^2_{\pi,\Delta\theta} (Z)]$. Denote $\nu= \frac{\|\Delta\theta\|}{\sqrt p}\mu_n^T(AA^T)^{-1}$, $\nu\in\bR^{n-1}$. We have
\begin{align*}
\E_0 [L^2_{\pi,\Delta\theta} (Z)] &=e^{-\rho^2} {{d\choose p}}^{-2} \!\!\!\!\! \sum_{m_1,m_2\in\cM(d,p)} 
\E_0 \Bigl[ \prod_{j\in m_1}\cosh\Bigl(\nu Z_j\Bigr)  \prod_{j\in m_2}\cosh\Bigl(\nu Z_j\Bigr)\Bigr] \\
&=e^{-\rho^2} {{d\choose p}}^{-2} \!\!\!\!\! \sum_{m_1,m_2\in\cM(d,p)} 
\E_0 \Bigl[ \prod_{j\in m_1\cap m_2}\cosh^2\Bigl(\nu Z_j\Bigr)  \prod_{j\in m_1\triangle m_2}\cosh\Bigl(\nu Z_j\Bigr)\Bigr] \\
&=e^{-\rho^2} {{d\choose p}}^{-2}\!\!\!\!\!  \sum_{m_1,m_2\in\cM(d,p)} 
\E_0 \Bigl[\Bigl(\cosh^2 (\nu Z_j)\Bigr)^{|m_1\cap m_2|} \Bigl(\cosh\nu Z_j\Bigr)^{|m_1\triangle m_2|}\Bigr].
\end{align*}
The following relations hold for $\xi\sim\cN(0,1)$, $\lambda,\theta\in\bR$,
$$
\E\cosh(\lambda(\theta+\xi))=e^{\lambda^2/2+\lambda\theta},\quad
\E\cosh^2(\lambda(\theta+\xi))=e^{\lambda^2+\lambda\theta}\cosh(\lambda^2+\lambda\theta).
$$
Noting that 
$$
|m_1\cap m_2|+|m_1 \triangle m_2|/2=p
$$
and taking into account that $\tilde \xi_j=A\xi_j$ and, consequently, $\nu Z_j\sim \cN(0, \rho^2/p)$ under $H_0$, we obtain
$$
\E_0 [L^2_{\pi,\Delta\theta} (Z)] =  {{d\choose p}}^{-2} \sum_{m_1,m_2\in\cM(d,p)} \Bigl[\cosh \Bigl(\frac{\|\Delta\theta\|^2}p nh(\tau)\Bigr)\Bigr]^{|m_1\cap m_2|}
$$
We need to study the behaviour of this sum for different values of $p$. 

Since the number of subsets that overlap by $j$ elements equals
$$
|(m_1,m_2)\in \cM(d,p)^2:\ m_1\cap m_2=j|={d\choose p}{p\choose j}{d-p\choose p-j}
$$
this expectation writes as 
$$
\E_0 [L^2_{\pi,\Delta\theta} (Z)] =  \sum_{j=0}^p \Bigl[\cosh \Bigl(\frac{\|\Delta\theta\|^2}p nh(\tau)\Bigr)\Bigr]^{j} p_{j,p,d}
$$
where $p_{j,p,d}$ is the hypergeometric distribution,
$$
p_{j,p,d}=\frac{{p\choose j}{d-p\choose p-j}}{{d\choose p}}.
$$
We now need to investigate the behaviour of
\begin{equation}\label{HGG}
\E_0 [L^2_{\pi,\Delta\theta} (Z)] =  \E  \Bigl[\cosh \Bigl(\frac{\|\Delta\theta\|^2}p nh(\tau)\Bigr)\Bigr]^{W},
\end{equation}
where $W$ is a hypergeometric $H(d,p,p)$ random variable. 

We recall here an inequality due to~\cite{Hoeffding:1963}; see also Appendix A.9 in~\cite{Devroye:Gyorfi:Lugosi:1996} for a recent reference.
If $W$ is a hypergeometric  with parameters $p,p,d$ and $Y$ is a binomial r.v. with parameters $p$ and $p/d$, then for any convex function $f$
\begin{equation}\label{Hoefding}
\E f(W)\le \E f(Y).
\end{equation}
Therefore, (\ref{HGG}) is bounded by
$$
\E_0 [L^2_{\pi,\Delta\theta} (Z)] \le   \E  \Bigl[\cosh \Bigl(\frac{\|\Delta\theta\|^2}p nh(\tau)\Bigr)\Bigr]^{Y}.
$$

Second, let us consider the ``Bernoulli prior'', that is the prior defined in~(ii) of this lemma. We follow a similar reasoning  as before,
$$
L^2_{\pi,\Delta\theta} (Z)
= {{d\choose p}}^{-1} e^{-\rho^2/2} \sum_{m\in\cM(d,p)} \prod_{j\in m}\exp\Bigl(-\frac{\|\Delta\theta\|}{\sqrt p}\mu_n^T(AA^T)^{-1} Z_j\Bigr)
$$
Now we calculate $\E_0 [L^2_{\pi,\Delta\theta} (Z)]$ and bound it from above using the same inequality~(\ref{Hoefding}),
\begin{align*}
\E_0 [L^2_{\pi,\Delta\theta} (Z)] &=e^{-\rho^2} {{d\choose p}}^{-2}\!\!\!\!\! \sum_{m_1,m_2\in\cM(d,p)} 
\E_0 \Bigl[ \prod_{j\in m_1} e^{-\nu Z_j} \prod_{j\in m_2} e^{-\nu Z_j}\Bigr] \\
&=e^{-\rho^2} {{d\choose p}}^{-2}\!\!\!\!\!  \sum_{m_1,m_2\in\cM(d,p)} 
\E_0 \Bigl[ \prod_{j\in m_1\cap m_2}e^{-2\nu Z_j}\prod_{j\in m_1\triangle m_2}e^{-\nu Z_j}\Bigr] \\
&=e^{-\rho^2} {{d\choose p}}^{-2} \!\!\!\!\! \sum_{m_1,m_2\in\cM(d,p)} 
\E_0 \Bigl[e^{-2\nu Z_j |m_1\cap m_2|} e^{-\nu Z_j |m_1\triangle m_2|}\Bigr]\\
&= {{d\choose p}}^{-2} \!\!\!\!\! \sum_{m_1,m_2\in\cM(d,p)} \exp \Bigl(\frac{\|\Delta\theta\|^2}p nh(\tau)|m_1\cap m_2|\Bigr)\\
& \le \E  \Bigl[\exp \Bigl(\frac{\|\Delta\theta\|^2}p nh(\tau)\Bigr)\Bigr]^{Y}.
\end{align*}
The lemma is proved. \hfill $\blacksquare$

The rest of the proof of Theorem~\ref{low_thm} follows from Lemmas~\ref{lowb1} and~\ref{lowb2}. 
\begin{lemma}\label{lowb1}
Let $\beta\in (0,1/2]$ and condition~(\ref{det_bound_lin}) in Theorem~\ref{low_thm} be satisfied. Then relation~(\ref{L2conv}) holds.
\end{lemma} 
{\bf Proof}. 
Let us put the prior $\P_{m,r}$ on the vector $\Delta\theta\in \Theta_{p,\tau}[r]\subseteq V_p^d]$ as defined in~(i) of Lemma~\ref{lemLRHG}. Then for a binomial r.v.\ $Y$ from Lemma~\ref{lemLRHG} we have 
\begin{align*}
\E_0 [L^2_{\pi,r} (Z)] &\le   \E  \Bigl[\cosh \Bigl(\frac{r^2}p nh(\tau)\Bigr)\Bigr]^{X}\\
&=\biggl(1+\frac pd\Bigl[\cosh \Bigl(\frac{r^2}p nh(\tau)\Bigr)-1\Bigr]\biggr)^p\\
&\le \exp\biggl(\frac{p^2}d\Bigl[\cosh \Bigl(\frac{r^2}p nh(\tau)\Bigr)-1\Bigr] \biggr).
\end{align*}
Since $\beta\le 1/2$, we have $d^{1-\beta}\ge \sqrt d$. Since $p\asymp d^{1-\beta}$, from condition~(\ref{det_bound_lin}) of Theorem~\ref{low_thm} it follows that
$$
\frac{r^4 n^2 h^2(\tau)}{2p^2} \le \frac{r^4 n^2 h^2(\tau)}{2d}\to 0,\quad d\to\infty.
$$
Using this property and the facts that $\cosh x^2\le \exp(x^4/2)$ and $e^x-1\le x+\frac12 x^2 e^x$ for $0<x<1$, we have
\begin{align*}
\E_0 [L^2_{\pi,r} (Z)]&\le \exp\biggl(\frac{p^2}d\Bigl[\exp \Bigl(\frac{r^4n^2 h^2(\tau)}{2p^2} \Bigr)-1\Bigr] \biggr)\\
&\le \exp\biggl(\frac{p^2}d\Bigl[\frac{r^4n^2 h^2(\tau)}{2p^2} +\frac12 \frac {r^8n^4h^4(\tau)}{4p^4}\exp\Bigl(\frac{r^4n^2h^2(\tau)}{2p^2}\Bigr)\Bigr] \biggr)\\
&\le \exp\biggl(\frac{r^4n^2 h^2(\tau)}{2d} +\frac {r^8n^4h^4(\tau)}{8d ^2}\exp\Bigl(\frac{r^4n^2h^2(\tau)}{2d}\Bigr)\biggr)\to 1,\quad d\to\infty
\end{align*}
and the lemma is proved. \hfill $\blacksquare$

\begin{lemma}\label{lowb2}
Let $\beta\in(1/2,1)$ and condition~(\ref{det_bound_scan}) of Theorem~\ref{low_thm} be satisfied.
Then relation~(\ref{L2conv}) holds.
\end{lemma}
{\bf Proof}. Let us put the Bernoulli prior $\P_{m,r}$ on the vector $\Delta\theta\in \Theta_{p,\tau}[r]\subseteq V_p^d]$ as defined in~(ii) of Lemma~\ref{lemLRHG}. 
Since $\beta>1/2$, $d^{1-\beta}<\sqrt d$.  Applying Lemma~\ref{lemLRHG} and using condition~(\ref{det_bound_scan}) of Theorem~\ref{low_thm}, we obtain for $p\asymp d^{1-\beta}$ 
\begin{align*}
\E_0 [L^2_{\pi,\Delta\theta} (Z)] &\le \exp\biggl(\frac{p^2}d\Bigl[\exp \Bigl(\frac{r^2 nh(\tau)}p\Bigr)-1\Bigr] \biggr)\\
&< \exp\biggl(\frac{p^2}d \left(\frac {d}p\right)^{2-1/\beta}\biggr)=\exp\biggl(\frac{p^{1/\beta}}{d^{1/\beta-1}} \biggr)\to 1,\quad d\to\infty.
 \end{align*}
The lemma is proved.\hfill $\blacksquare$

\section{Auxiliary results}\label{app}

The following result from~\cite{Spokoiny:2013}  on the deviation bounds for quadratic forms is used throughout the paper.

\begin{lemma}\label{Spok}
Let $\xi$ be a standard normal vector in $\bR^p$. Then for any $u>0$ it holds
$$
\P(\|\xi\|^2>p+u) \le \exp\{-\frac p2\varphi(u/p)\}
$$
with $\varphi(t)=t-\log(1+t)$. For any $x$,
$$
\P(\|\xi\|^2>p+\varphi^{-1}(2x/p))\le \exp(-x).
$$
This particularly yields with $\varkappa=6.6>2/(1-\log 2)$
\begin{equation}\label{ineq_spok}
\P(\|\xi\|^2>p+\sqrt{\varkappa xp}\vee (\varkappa x))\le \exp(-x).
\end{equation}
\end{lemma}

Below we give the proof of formula~(\ref{LLZ})  in Remark~1 regarding the log-likelihood.
\begin{lemma}\label{LL_lemma}
For observations~(\ref{channels}) the marginal log-likelihood of $\tau$ and $m$ is given by relation~(\ref{LLZ}).
\end{lemma}

{\bf Proof}. 
Denote  $X_i=(X_i^1,\dots,X_i^d)$, $i=1,\dots,n$ and $X^j=(X_1^j,\dots,X_n^j)$, $j=1,\dots,d$. We can rewrite our data~(\ref{channels}) in the following way,
\begin{equation}\label{data}
X_i^j=\theta^j+\delta_j\Delta\theta_\tau^j\1\{i>\tau\}+\xi_i^j,\quad i=1,\dots,n,\quad j=1,\dots,d,
\end{equation}
where $\delta_j$ are binary variables taking the values 1 or 0 depending on whether there is a change in the component $j$ or not such that 
$$
\delta_j=\begin{cases} 
1,& \mbox{there is a change, $j\in m$}\\
0,& \mbox{there is no change, $j\notin m$}.
\end{cases}
$$
Let $\P^j_{\theta^j,\Delta\theta_\tau^j\delta_j}$ be the measure corresponding to the case of a change in the data with respective means $\theta^j$ and $\theta^j+\Delta\theta_\tau^j$ and $\P_0$ be the standard Gaussian measure. Then the negative log-likelihood ratio of estimating $(\theta^j,\Delta\theta_\tau^j\delta_j,\tau)$ from the sequence $X^j$ is given by
\begin{equation}\label{LL}
-\log \frac{d\P^j_{\theta^j,\Delta\theta_\tau^j\delta_j}}{d\P_0} (X^j)=\frac1{2}\sum_{i=1}^\tau (X_i^j-\theta_j)^2+\frac1{2}\sum_{i=\tau+1}^n (X_i^j-\delta_j\Delta\theta_j-\theta_j)^2.
\end{equation}
Estimating the values of $\theta_j$ and $\Delta\theta_j$, 
$$
\widehat\theta_j=(\tau+(1-\delta_j)(n-\tau))^{-1}\sum_{i=1}^\tau X_i^j,\quad
\widehat\theta_j+\widehat{\Delta\theta_j}=(n-\tau)^{-1}\delta_j\sum\limits_{i=\tau+1}^n X_i^j,
$$
and then substituting them into~(\ref{LL}) we obtain 
\begin{align*}
-\log\frac{d\P^j_{\delta_j,\tau}}{d\P_0}=-\frac{\delta_j}{2}[Z_n^j(\tau)]^2,
\end{align*}
where 
$$
[Z_n^j(\tau)]^2=\frac1\tau \biggl(\sum_{i=1}^\tau X_i^j\biggr)^2+\frac1{n-\tau} \biggl(\sum_{i=\tau+1}^n X_i^j\biggr)^2-\frac1n\biggl(\sum_{i=1}^n X_i^j\biggr)^2
$$
are as defined in~(\ref{Zj}). 
Next, from independence of the vectors $X^j$ we obtain 
$$
\log \frac{d\P_{m,\tau}}{d\P_0}(X) =\sum_{j=1}^d \log \frac{d\P^j_{\delta_j,\tau}}{d\P^j_0}(X^j)= \frac1{2} \sum_{j=1}^d \delta_j [Z_n^j(\tau)]^2
=\frac12 \sum_{j\in m} [Z_n^j(\tau)]^2
$$
where $m$ is a set of components with a change. \hfill $\blacksquare$

\section*{Acknowledgements}
The authors would like to thank Alexander Goldenshluger, Yuri Golubev, Anatoli Juditsky, Axel Munk, and David Siegmund for fruitful discussions.

\end{document}